# OPTIMAL FLOW THROUGH THE DISORDERED LATTICE[1]

By David Aldous

*University of California at Berkeley*

Consider routing traffic on the $N \times N$ torus, simultaneously between all source-destination pairs, to minimize the cost $\sum_e c(e) f^2(e)$, where $f(e)$ is the volume of flow across edge $e$ and the $c(e)$ form an i.i.d. random environment. We prove existence of a rescaled $N \to \infty$ limit constant for minimum cost, by comparison with an appropriate analogous problem about minimum-cost flows across a $M \times M$ subsquare of the lattice.

**1. Introduction.** In highly abstracted models of transportation or communication (e.g., roads, Internet) one is required (simultaneously for all source-destination pairs) to route a certain "volume of flow" from source to destination, and one seeks to minimize some notion of cost subject to some constraints (e.g., edge-capacities). In contrast to queueing theory, we shall regard flows as deterministic, but networks as random. A survey of such problems under various different models of random networks will be given elsewhere. In this paper we focus on the particular setting of *disordered* $\mathbb{Z}^2$ and, after initial general discussion, study a specific model of costs.

The mathematical structure of the lattice $\mathbb{Z}^2$ with an i.i.d. random environment $(c(e) : e$ an edge of $\mathbb{Z}^2)$ arises in many settings, such as first passage percolation [10, 14], disordered Ising models [4] and random walk in random environment [20]. Conceptually, the point of such *disordered lattice* models is to represent local spatial inhomogeneity. To study flow problems, consider a $N \times N$ square within $\mathbb{Z}^2$. To see how quantities scale suppose, for each pair $(v_1, v_2)$ of vertices within the square, we send a flow of volume $a_N$ from $v_1$ to $v_2$ by spreading it evenly amongst minimal-length paths. Then the total volume of flow along a typical edge will be of order

$$a_N \times N^4 \times N/N^2 = N^3 a_N,$$

Received December 2005; revised March 2006.
[1]Supported in part by NSF Grant DMS-02-03062.

where number of source-destination pairs $= N^4$; path length is of order $N$; number of edges is of order $N^2$. We therefore scale by writing $a_N = \rho N^{-3}$, and we interpret $\rho$ as normalized traffic intensity. This normalization convention means that, in the optimal flow problems, we shall consider, flow volume along individual edges should be of order 1. To formulate one problem, suppose we impose capacity constraints: the total flow $f(e)$ along each edge $e$ cannot exceed $\mathrm{cap}(e)$, where the $(\mathrm{cap}(e):e$ an edge of $\mathbb{Z}^2)$ are i.i.d. with support bounded away from 0 and infinity. It seems intuitively very plausible that there should be a critical value $\rho_0$ such that

$$\mathbb{P}(\text{feasible flow exists}) \to 1 \qquad \text{as } N \to \infty, \ \rho < \rho_0$$

$$\mathbb{P}(\text{feasible flow exists}) \to 0 \qquad \text{as } N \to \infty, \ \rho > \rho_0.$$

Moreover, it seems plausible that, if we also have i.i.d. costs $(\mathrm{cost}(e):e$ an edge of $\mathbb{Z}^2)$ representing cost per unit volume of flow across edge $e$, and if (when a feasible flow exists) we consider the cost $C_N$ of a minimum-cost feasible flow, then we should have a limit function

$$\lim_N N^{-2} \mathbb{E} C_N = g(\rho), \qquad \rho < \rho_0.$$

These are the kinds of results we would like to prove, but in fact (to avoid a specific technical issue discussed in Section 4.1) we replace "hard constraints" $\mathrm{cap}(e)$ on the flow through an edge by "soft constraints" as follows. There is an i.i.d. environment $(c(e):e$ an edge of $\mathbb{Z}^2)$ whose interpretation is now the cost of a total volume $f(e)$ of flow across edge $e$ equals $c(e)f^2(e)$.

So the total cost of a flow is defined as

(1) $$\sum_e c(e) f^2(e)$$

summed over the edges $e$ of the $N \times N$ square. Take normalized traffic intensity $\rho = 1$, and consider the cost $C_N$ of a minimum-cost flow ($C_N$ is a random variable, the randomness arising only from the random environment). Because flow volume across a typical edge should be of order 1, we expect a limit constant

$$\lim_N N^{-2} \mathbb{E} C_N = \gamma,$$

depending only on the distribution of $c(e)$. This is essentially what our main result (Theorem 1) establishes, with two further provisos. To exploit global stationarity, we work on the $N \times N$ torus instead of the square; and for technical reasons (Section 4.1) we impose a constant bound $B$ on volume of flow across edges.

Conceptually speaking, we are modeling *congestion* in the following sense. In a road network, the cost to users is the time taken to traverse a route. In the absence of congestion the cost per user does not depend on the flow volume, so the total cost (all users combined) associated with an edge $e$ is linear in $f(e)$. If congestion causes speeds to decrease, it increases costs to all users of an edge, so the total cost is superlinear in $f(e)$. The choice of $f^2(e)$ in (1) is just an arbitrary but mathematically convenient choice of superlinear function.

Merely proving existence of a limit, as this paper does, is of course a rather modest achievement, so let us attempt several justifications:

(i) Experience with other "disordered lattice" problems such as first passage percolation suggests it is very difficult to obtain explicit formulas for limit constants.

(ii) A natural way to describe network behavior is via functions $g_N(\rho)$ describing how some quantity varies with traffic intensity $\rho$. In the context of probability models for $N$-vertex networks, one would like to ensure existence of suitably renormalized $N \to \infty$ limit functions $g(\rho)$, and this paper is part of a project exploring methodology for such results.

(iii) The standard technique for proving existence of limits in "spatial optimization" problems is subadditivity, which is readily applicable to, for example, traveling salesman type problems [18, 19], but which does not apply easily to our "flow" problems, because there is no simple way to relate the optimal flow on a $2N \times 2N$ square to the optimal flows on the four $N \times N$ subsquares. Instead we use the *local weak convergence* methodology explained in Section 2.1. Our interest in this methodology arose from its use in optimization problems over locally tree-like random networks, where it becomes a reformulation of (special settings of) the *cavity method* [16] of statistical physics and allows explicit albeit nonrigorous calculation of limit constants.

*Plan of paper.* Section 2 sets out notation, states the results carefully, and outlines the proof. The proofs themselves comprise Sections 3 and 4. We conclude in Section 5 with a wide-ranging discussion—about optimal flows in general, about local weak convergence methodology, about details of our particular model and about related work.

**2. Results and notation.** Section 2.1 gives a rough verbal overview of methodology. We set up general notation in Section 2.2, enabling us to state our convergence result (Theorem 1) in Section 2.3. In Section 2.4 we set up notation specifically for dealing with flows across squares in the lattice, which permits us to state Theorem 2 encapsulating the methodological idea of relating flows within the $N \times N$ torus to flows across the $M \times M$ square. Section 2.5 then provides a more detailed outline of the proofs to follow.

2.1. *Methodology.* For $M < N$, it is obvious that an i.i.d. environment in a $M \times M$ subsquare of the $N \times N$ torus is distributed as an i.i.d. environment in a $M \times M$ subsquare of the lattice $\mathbb{Z}^2$. Our methodology rests on the less obvious idea that the $N \to \infty$ limit of the global cost $C_N$ on the torus (which involves route-lengths of order $N$), suitably normalized, can be identified with the $M \to \infty$ limit of a certain quantity $c_{M,B}$ defined in terms of flows across the $M \times M$ square. We say *across* (from boundary to boundary) because as $N \to \infty$ the volume of flow originating at a vertex becomes negligible compared to the volume passing through the vertex. The argument has three parts. For a flow across a square, we call the joint distribution of entrance and exit points the *transportation measure*.

*Part* 1. We consider $N \to \infty$ limits of the empirical distribution, over all $M \times M$ subsquares, of the transportation measure for the optimal flow, jointly with the environment.

*Part* 2. We use a concentration of measure argument to show that, given the empirical distribution of the transportation measure over a $M \times M$ square, for most realizations of the $(c(e))$ within the square the cost of optimal flow across the square is close to the mean cost. This allows us to take expectation over environments. Then convexity of the cost function (1) allows us to replace the empirical transportation measure by its mean measure.

*Part* 3. The argument above leads to a definition (10) of $c_{M,B}$ which provides *lower* bounds on the $N \to \infty$ limit optimal cost. To get the corresponding upper bound, we need to construct flows on the $N \times N$ torus. Partition the torus into $M \times M$ squares. The definition of $c_{M,B}$ involves some particular transportation measure $Q_0$ across $M \times M$ squares. We construct flows by first constructing a "skeleton" as a Markov chain which steps from one boundary (between two squares) point to another boundary point; this skeleton does not depend on the realization of the environment. Within each square, use a minimum cost flow (which does depend on realization of environment) consistent with the transportation measure.

2.2. *General notation.* Consider a graph $G = (\mathcal{V}, \mathcal{E})$, which for our purposes will usually be either the discrete $N \times N$ torus $[0, N-1]^2$ (which we call $\mathcal{T}_N^2$) or the extended $M \times M$ square described in Section 2.4. A finite oriented path $\sigma = (v_0, v_1, \ldots, v_k)$ has a source or *entrance* $\text{ent}(\sigma) = v_0$ and a destination or *exit* $\text{exi}(\sigma) = v_k$. For an edge $e$ of $G$, write $n(\sigma, e)$ for the number of times that $e$ occurs (in either direction) in $\sigma$. Almost always we implicitly deal with self-avoiding paths, for which $n(\sigma, e) = 0$ or 1. A *path-flow* is a measure (by which we always mean a *nonnegative finite* measure) $\mu$ on the space $\Sigma$ of finite paths. The map

$$\sigma \to (\text{ent}(\sigma), \text{exi}(\sigma))$$

from $\Sigma$ to $\mathcal{V} \times \mathcal{V}$ pushes $\mu$ forward to a measure on $\mathcal{V} \times \mathcal{V}$ which we call the *transportation measure* $\text{tra}(\mu)$ by analogy with the classical mass transportation problem [17]. Call the edge-indexed collection

$$\mathbf{f} = (f(e), e \in \mathcal{E}) = \left( \int n(\pi, e) \mu(d\pi), \ e \in \mathcal{E} \right)$$

the *flow-volume* $\text{flo}(\mu)$ associated with $\mu$. Thus, if $\text{tra}(\mu)$ concentrated on a single element $(v, v')$, then $\text{flo}(\mu)$ would be a flow from source $v$ to destination $v'$ in the elementary sense of the max-flow min-cut theorem. But we use "flow" in the sense of *multicommodity flow*, and in informal discussion we do not distinguish carefully between path-flows and their associated flow-volumes.

Consider a collection $\mathbf{c} = (c(e) : e$ an edge of $\mathcal{E})$ of nonnegative real numbers. Call $c(e)$ a *cost-factor* and call $\mathbf{c}$ an *environment*. Given a flow-volume $\mathbf{f}$ and an environment $\mathbf{c}$, define the cost of that flow in that environment to be

$$\text{cost}(\mathbf{f}, \mathbf{c}) = \sum_{e \in \mathcal{E}} c(e) f^2(e). \tag{2}$$

To make our probability model, let the environment $\mathbf{c}$ be chosen i.i.d. from some distribution $\kappa$ with bounded support: for some $0 < c^* < \infty$,

$$\kappa[0, c^*] = 1; \qquad (c(e) : e \text{ an edge of } \mathcal{E}) \text{ are i.i.d. } (\kappa). \tag{3}$$

In what follows, $\mathbb{E}$ and $\mathbb{P}$ denote expectation and probability with respect to the random environment.

Call a measure $\theta$ on a finite set $B$ *constant* if $\theta(b)$ is constant for all $b \in B$; we reserve *uniform* for the case of a probability measure.

Write $M | N$ for "$N$ is a multiple of $M$."

2.3. *The main theorem.* Now take the graph to be the discrete torus $\mathcal{T}_N^2 = [0, N-1]^2$, and where helpful append a subscript $\cdot_{(N)}$ to the notation above, so that, for instance, $\Sigma_{(N)}$ is the set of finite paths in $\mathcal{T}_N^2$. As explained informally in the Introduction, we want to consider flows with volume $N^{-3}$ between every source-destination pair. In the terminology above, define a *standardized global flow* $\mathbf{f}$ to be the flow-volume of some path-flow $\mu$ for which $\text{tra}(\mu)$ is the constant measure

$$\text{tra}(\mu)(v, w) = N^{-3}, \qquad (v, w) \in \mathcal{T}_N^2 \times \mathcal{T}_N^2.$$

Fix $B > \frac{1}{4}$ for the rest of the paper. Then (Lemma 6) there exist standardized global flows satisfying the capacity constraint

$$\max_e f(e) \leq B. \tag{4}$$

For a fixed environment **c**, consider the minimum cost subject also to this capacity constraint:

$$\text{cost}_{(N)}(\mathbf{c}, B) := \inf\{\text{cost}(\mathbf{f}, \mathbf{c}) : \mathbf{f} \text{ a standardized global flow satisfying (4)}\}.$$

Now make the environment random as at (3).

THEOREM 1. *There exists a constant $\gamma(\kappa, B)$ such that*

$$\lim_{N \to \infty} N^{-2} \mathbb{E} \, \text{cost}_{(N)}(\mathbf{c}, B) = \gamma(\kappa, B).$$

The corresponding SLLN then holds by a concentration inequality; see Section 5.

2.4. *Flows across subsquares*: *notation.* We now develop notation for use in the proofs.

Figure 1 (middle right) shows the $M \times M$ square $[0, M-1]^2$ (regarded as a subgraph of the lattice $\mathbb{Z}^2$) together with half of each edge from this square to its complement. Create artificial *boundary points* $b$ at ends of these half-edges; each such $b$ is of the form $(-\frac{1}{2}, j)$ or $(M - \frac{1}{2}, j)$ or $(i, -\frac{1}{2})$ or $(i, M - \frac{1}{2})$. Write $\mathcal{B}_M^o$ for the set of all $M^2$ internal vertices. Write $\text{Bou}_M$ for the set of $4M$ boundary points, and write $\mathcal{B}_M = \mathcal{B}_M^o \cup \text{Bou}_M$. Write $\mathcal{E}_M$ for the edge-set (internal edges and boundary half-edges). Call $(\mathcal{B}_M, \mathcal{E}_M)$ the *extended $M \times M$ square*. We use the general notation of Section 2.2 for the graph $(\mathcal{B}_M, \mathcal{E}_M)$, writing a subscript $\cdot_M$ where helpful. So $\Sigma_M$ denotes the set of finite paths in this graph, and a path-flow $\mu$ has a transportation measure $\text{tra}(\mu)$ on $\mathcal{B}_M \times \mathcal{B}_M$ and a flow-volume $\mathbf{f} = \text{flo}(\mu)$ on the space $\mathcal{F}_M \subset [0, \infty)^{\mathcal{E}_M}$ of possible flow-volumes. Write $\mathbf{c} = (c(e), e \in \mathcal{E}_M)$ for an environment of cost-factors on the edges of the extended $M \times M$ square; write $\mathcal{C}_M$ for the set of all such environments $\mathbf{c}$. For a flow-volume $\mathbf{f}$ and an environment $\mathbf{c}$, we write

$$\text{cost}_M(\mathbf{f}, \mathbf{c}) = \sum_{e \in \mathcal{E}_M} c(e) f^2(e). \tag{5}$$

The preceding just copied general notation; we now come to new definitions specific to the extended $M \times M$ square. Partition the boundary points $\text{Bou}_M$ into four subsets

$$\text{Bou}_M = \text{Bou}_M^{\text{top}} \cup \text{Bou}_M^{\text{bottom}} \cup \text{Bou}_M^{\text{left}} \cup \text{Bou}_M^{\text{right}}$$

in the way implied by the notation (Figure 1, bottom right). For $b \in \text{Bou}_M$, let $b^{\text{reflect}} \in \text{Bou}_M$ be the boundary point obtained by reflecting $b$ top-to-bottom or left-to-right; in particular, $b^{\text{reflect}} = b \mod (M, M)$ (Figure 1, bottom right).

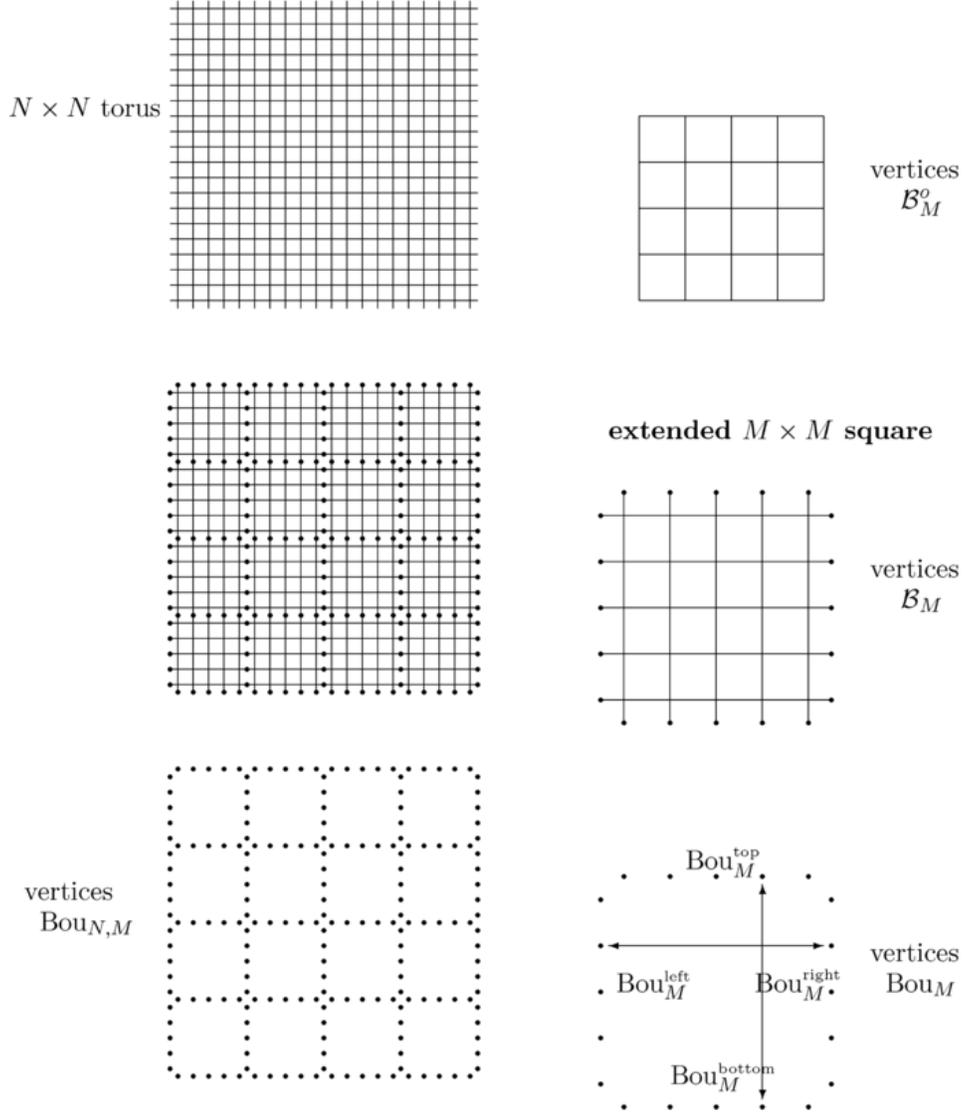

FIG. 1. Top left: *the $N \times N$ torus (here $N = 20$), in which the edges wrap from left side to right side and from top to bottom.* Middle left: *the natural partition of the $N \times N$ torus into $M \times M$ subsquares. Here $M = 5$.* Bottom left: *the boundary points $\mathrm{Bou}_{N,M}$ for the partition.* Top right: *the "ordinary" $M \times M$ square $[0, M-1]^2$.* Middle right: *the $M \times M$ square with its boundary points. This is what we call "the extended $M \times M$ square" in the text.* Bottom right: *The set $\mathrm{Bou}_M$ of boundary points of the extended $M \times M$ square, and the map $b \to b^{\mathrm{reflect}}$.*

Given a transportation measure $Q$ on $\text{Bou}_M \times \text{Bou}_M$, write $Q_{\text{ent}}$ and $Q_{\text{exi}}$ for its marginals on $\text{Bou}_M$ (later we use the same notation for transportation measures on other grids). Define $\mathcal{Q}_M$ to be the set of transportation measures $Q$ such that

(6) $\qquad$ the push-forward of $Q_{\text{exi}}$ under $b \to b^{\text{reflect}}$ equals $Q_{\text{ent}}$.

For $Q \in \mathcal{Q}_M$, define

$$\text{drift}(Q) = \frac{1}{M^2} \sum_{b_1} \sum_{b_2} (b_2 - b_1) Q(b_1, b_2). \tag{7}$$

So $\text{drift}(Q)$ is a point in $\mathbb{R}^2$. Taking this point modulo $(1,1)$ gives a point in the continuous torus $\mathcal{T}^2 = [0,1)^2$:

$$\text{drift}_1(Q) := \text{drift}(Q) \bmod (1,1).$$

Finally, define $Q_0$ to be *isotropic* if it has a representation as a mixture

$$Q_0 = \int_{\mathcal{Q}_M} Q \, \psi(dQ), \tag{8}$$

where $\psi$ is a probability distribution on $\mathcal{Q}_M$ whose push-forward under the map $Q \to \text{drift}_1(Q)$ is the uniform probability distribution on $[0,1)^2$. Write $\mathcal{Q}_M^{\text{iso}}$ for the set of isotropic transportation measures. Because property (6) is preserved under mixtures, $\mathcal{Q}_M^{\text{iso}} \subset \mathcal{Q}_M$.

Given a transportation measure $Q$ on $\text{Bou}_M \times \text{Bou}_M$, and an environment $\mathbf{c}$, define

$$\text{cost}_{M,B}(Q, \mathbf{c}) = \inf_{\mu} \left\{ \text{cost}_M(\text{flo}(\mu), \mathbf{c}) : \text{tra}(\mu) = Q, \max_e \text{flo}(\mu)(e) \leq B \right\}. \tag{9}$$

We now arrive at the central definition:

$$c_{M,B} := \inf \{ \mathbb{E} \, \text{cost}_{M,B}(Q, \mathbf{c}) : Q \in \mathcal{Q}_M^{\text{iso}} \}. \tag{10}$$

Here $\mathbb{E}$ is expectation w.r.t. the random environment, with the following convention. Cost-factors $c(e)$ are i.i.d. ($\kappa$) as $e$ runs over internal edges and over edges to boundary points on the bottom and left sides; $c(e) = 0$ on other edges to boundary points. Obviously, this convention is designed so that when the $N \times N$ torus is partitioned into $M \times M$ squares (Figure 1, middle left), each edge of the torus is assigned to a unique square. In bounding costs we may take the number of edges $|\mathcal{E}_M|$ as $2M^2$.

Here is the promised elaboration of Theorem 1.

THEOREM 2. *For $B > 1/4$, there is a limit constant*

$$\gamma(\kappa, B) = \lim_{M \to \infty} M^{-2} c_{M,B}.$$

*And*

$$\lim_{N \to \infty} N^{-2} \mathbb{E} \, \text{cost}_{(N)}(\mathbf{c}, B) = \gamma(\kappa, B).$$

2.5. *Outline of proof.* Here we expand the rough description of methodology (Section 2.1) into a more detailed outline of the proof. This outline will ignore the constraint $f(e) \leq B$, which enters only at a technical level, sometimes as a convenience and sometimes as an inconvenience.

Getting an upper bound on the limit $\lim_{N\to\infty} N^{-2} \mathbb{E} \operatorname{cost}_{(N)}(\mathbf{c}, B)$ requires a construction of a standardized global flow. The definition (10) of $c_{M,B}$ involves some particular transportation measure $Q \in \mathcal{Q}_M^{\text{iso}}$ attaining the infimum. This $Q$ is a positive matrix on the boundary vertices of a standard extended $M \times M$ square, but can be used in a natural way to define a Markov chain on the boundary points $\operatorname{Bou}_{N,M}$ of the partition of the $N \times N$ torus into $M \times M$ squares: see Figure 2. Condition (6) is the essential "compatibility" condition making this construction work. Running the chain for a suitable number of steps with a suitable starting distribution on $\operatorname{Bou}_{N,M}$ gives a joint distribution of starting and ending points which is approximately uniform on $\operatorname{Bou}_{N,M} \times \operatorname{Bou}_{N,M}$. The definition (8) of *isotropic* is what gives this "approximately uniform" property (Proposition 4). Making the walk into a standardized global flow requires appending length $o(N)$ path-segments at the ends to make the transportation measure become exactly uniform, and expanding a step between two points on the boundary of an $M \times M$ square into a path across the square. But the construction makes the transportation measure across a square be the original $Q$, so to "expand a step" we simply use the routing attaining the minimum in the definition (9) of $\operatorname{cost}_{M,B}(Q, \mathbf{c})$; see Proposition 5. This completes the construction of a standardized global flow. The cost associated with the end $o(N)$-length segments is negligible, and the cost associated with flow across a typical $M \times M$ square in the torus is $c_{M,B}$, leading to the desired upper bound (Proposition 17).

The lower bound uses more abstract methods. Consider the optimal flow within the $N \times N$ torus, then consider this flow within a randomly-positioned $M \times M$ square, and then consider subsequential $N \to \infty$ weak limits for fixed $M$. What can we say about such a limit distribution? At first sight we cannot say anything explicit, because we do not know anything explicit about the optimal flow on the $N \times N$ torus. But we can collect some properties, as follows:

(i) The environment $\mathbf{c} = (c(e))$ within the $M \times M$ square follows the original i.i.d. model (3).

(ii) Because the volume of flow within the torus originating or ending at a given vertex is $N^{-1}$, in the present $N \to \infty$ limit all the flow is *across* the $M \times M$ square from boundary to boundary.

(iii) There is some random transportation measure $\mathbf{Q}$ giving the joint distribution of entrance and exit points of flow. Note $\mathbf{Q}$ will be dependent on $\mathbf{c}$.

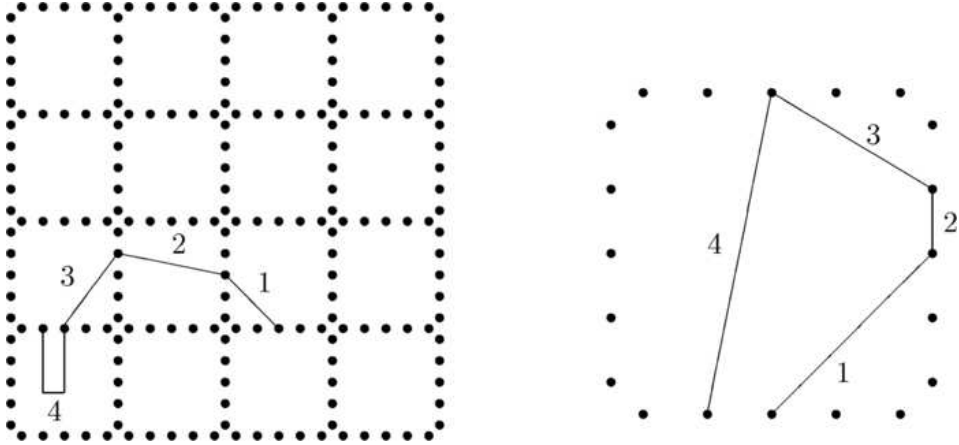

FIG. 2. Left side: 4 *steps* $(Y_0, Y_1, Y_2, Y_3, Y_4)$ *of a walk on the skeleton graph. Right side: the projection walk* $(\mathrm{pro}_{N,M}(Y_0), \mathrm{pro}_{N,M}(Y_1), \mathrm{pro}_{N,M}(Y_2), \mathrm{pro}_{N,M}(Y_3), \mathrm{pro}_{N,M}(Y_4))$ *on the boundary of the extended* $M \times M$ *square. The state of the projection walk is the relative position of entry to the next square for the skeleton walk.*

(iv) Because a global optimum is a local optimum, given the environment $\mathbf{c}$ and the transportation measure $\mathbf{Q}$, the flow-volumes $\mathbf{f} = (f(e))$ within the $M \times M$ square minimize the local cost

$$\mathrm{cost}_M(\mathbf{f}, \mathbf{c}) = \sum_{e \in \mathcal{E}_M} c(e) f^2(e)$$

subject to the given transportation measure.

(v) The fact that the optimal flow on the torus was some *standardized global flow* constrains the induced flows across a typical $M \times M$ square, and what it implies is (roughly, and in part) that $\mathbb{E}\mathbf{Q} \in \mathcal{Q}_M^{\mathrm{iso}}$.

(vi) The mean cost on the torus relates to the mean cost associated with our typical $M \times M$ square.

Lemma 20 formalizes these assertions in the finite $N$ context. Letting $N \to \infty$ gives a lower bound for $\lim_{N \to \infty} N^{-2} \mathbb{E} \mathrm{cost}_{(N)}(\mathbf{c}, B)$ involving

(11) $\qquad \mathbb{E} \inf\{\mathrm{cost}_M(\mathbf{f}, \mathbf{c});\ \mathbf{f}$ has transportation measure $\mathbf{Q}\}$.

Lemma 21 states this bound precisely. Of course, we do not know explicitly what is the distribution of $\mathbf{Q}$, so we need to lower bound the quantity (11) by appealing to constraints $\mathbf{Q}$ must satisfy. To do so, the key ingredients are the two very easy facts referenced below. Consider a nonrandom $Q$. Then by an easy concentration inequality (28), the random variable $\mathrm{cost}_{M,B}(Q, \mathbf{c})$ is close to its expected value $\mathbb{E} \mathrm{cost}_{M,B}(Q, \mathbf{c})$. So if $\mathbf{Q}$ were random but took values in a small set $(Q_j)$ of possible values, then the expectation in (11) would be approximately the appropriate weighted average of $\mathbb{E}\operatorname{cost}_{M,B}(Q_j, \mathbf{c})$. Because (Lemma 18) $Q \to \mathbb{E}\operatorname{cost}_{M,B}(Q, \mathbf{c})$ is convex, we could then lower bound the quantity (11) by $\mathbb{E}\operatorname{cost}_{M,B}(\bar{Q}, \mathbf{c})$ for $\bar{Q} = \mathbb{E}\mathbf{Q} \in \mathcal{Q}_M^{\text{iso}}$, which by definition is lower bounded by $c_{M,B}$, giving the desired lower bound in Theorem 2.

So the remaining issue is to show that the quantity (11) is not much changed by replacing $\mathbf{Q}$ by a "quantized" version taking values in some large finite set of "smoothed" transportation measures, which we define as follows. Set $M = KL$ and partition the boundary of the $M \times M$ square into blocks of length $K$. Require a smoothed transportation measure $Q'$ to have the property that $Q'(b_1, b_2)$ depends only on the blocks containing $b_1$ and $b_2$, and that $Q'(b_1, b_2)$ be an integer multiple of some small constant. Then it is enough to show that, in an arbitrary environment $\mathbf{c}$, the cost of the cheapest flow $\mathbf{f}$ with transportation measure $Q$ is not much less than the cost of the cheapest flow $\mathbf{f}'$ with smoothed transportation measure $Q'$. This scheme turns out to be technically complicated to implement (see Section 4.3), but a variant is carried though in Lemma 19.

**3. Constructing global flows from local flows.** In this section we show how to use flows across the extended $M \times M$ square to construct standardized global flows on the $N \times N$ torus. This leads to the upper bound of Proposition 17. The first stage of the argument leads to a "clean" construction of a global flow (Proposition 5) for which the source-destination distribution is *approximately* constant. We then need to "patch" to make the source-destination distribution be *exactly* constant; this is done using some of the elementary constructions collected in Section 3.2.

3.1. *The basic construction.* Fix $M|N$. Partition the $N \times N$ torus into $M \times M$ squares in the natural way (Figure 1, middle left). Insert boundary points midway along edges linking different squares. This creates a set $\operatorname{Bou}_{N,M}$ of $2N^2/M$ boundary points. Define the *skeleton graph* $(\operatorname{Bou}_{N,M}, \mathcal{E}_{N,M})$ to have an edge between each pair of boundary points which are boundary points of some common square $S$. Consider a measure $\nu$ on paths $(x_i)$ in the skeleton graph. For any square $S$ with boundary $\operatorname{Bou}_{(S)}$, and for each $(b_1, b_2) \in \operatorname{Bou}_{(S)} \times \operatorname{Bou}_{(S)}$, write

$$(12) \qquad \nu_{(S)}(b_1, b_2) = \sum_i \nu\{(x_0, x_1, x_2, \ldots) : (x_i, x_{i+1}) = (b_1, b_2)\}$$

for the mean number of times the path steps from $b_1$ to $b_2$. We often discuss a typical $M \times M$ square as if it were the standard $M \times M$ square by saying "up to translation."

We define some properties for a measure $\theta$ on $\operatorname{Bou}_{N,M} \times \operatorname{Bou}_{N,M}$ with marginals $\theta_{\text{ent}}$ and $\theta_{\text{exi}}$.

PROPERTIES 3.   (i) $\theta$ has total mass $N$.

(i) $\theta$ is invariant under translation of the torus by $(M, 0)$ or $(0, M)$.
(ii) $\theta_{\text{ent}} = \theta_{\text{exi}}$.

Recall that $\mathcal{Q}_M$ is a set of transportation measures on the boundary of the extended $M \times M$ square.

PROPOSITION 4.   *Given $Q \in \mathcal{Q}_M$ satisfying an irreducibility property (14), there is a measure $\theta_{N,M}$ on $\text{Bou}_{N,M} \times \text{Bou}_{N,M}$ satisfying Properties 3 and the following two properties:*

(i) *There is a measure $\nu$ on paths in the skeleton graph whose transportation measure $\text{tra}(\nu) = \theta_{N,M}$ and such that, for each square $S$, the measure $\nu_{(S)}$ at (12) equals $Q$ (up to translation).*

(ii) *Let $\widetilde{\theta}_{N,M}$ be the push-forward of $\theta_{N,M}$ under the map $(z^1, z^2) \to (\frac{1}{N} z^1, \frac{1}{N} z^2)$ from $\text{Bou}_{N,M} \times \text{Bou}_{N,M}$ to $\mathcal{T}^2 \times \mathcal{T}^2$ (where $\mathcal{T}^2$ is the continuous torus $[0, 1)^2$). Then as $N \to \infty$, there is weak convergence of $N^{-1} \widetilde{\theta}_{N,M}$ to the distribution of $(U, U + \text{drift}_1(Q))$, where $U$ has uniform probability distribution on $\mathcal{T}^2$.*

PROOF.   Use $Q$ to define a transition matrix $\bar{Q}$ on $\text{Bou}_M$ via

$$\bar{Q}(b, b') = Q(b, b')/Q_{\text{ent}}(b).$$

This is a (Markov) transition matrix because $Q_{\text{ent}}(b) = \sum_{b'} Q(b, b')$. The set $\text{Bou}_{N,M}$ of inter-square boundary points $b$ can be expanded to a set $\text{Bou}_{N,M}^+$ with elements

$$(b, \to) \text{ and } (b, \leftarrow); \quad \text{or} \quad (b, \uparrow) \text{ and } (b, \downarrow)$$

indicating a direction along an inter-square edge. Use the transition matrix $\bar{Q}$ to define a transition matrix $Q^1(\cdot, \cdot)$ on $\text{Bou}_{N,M}^+$ as follows. For a state, say, $(b_0, \to)$, the arrow points into some square $S_0 = [x_0, x_0 + M - 1] \times [y_0, y_0 + M - 1]$, and $b_0$ is of the form $b_0 = (x_0, y_0) + b_0^*$ for some $b_0^* \in \text{Bou}_M$. Any $b_1^* \in \text{Bou}_M$ determines a point $b_1 = (x_0, y_0) + b_1^* \in \text{Bou}_{N,M}$. Define

$$Q^1((b_0, \to), (b_1, \sharp)) = \bar{Q}(b_0^*, b_1^*),$$

where $\sharp$ is the direction of arrow pointing out of $S_0$ at $b_1$.

Using the defining property (6) for $Q \in \mathcal{Q}_M$, one can check that a stationary distribution $\pi$ on $\text{Bou}_{N,M}^+$ for $Q^1$ is

(13) $$\pi((x_0, y_0) + b_1^*, \sharp) = Q_{\text{ent}}(b_1^*)/q,$$

where as above $\sharp$ is the direction of arrow pointing out of $S_0$ at $b_1$, and where the normalization constant is $q = (\frac{N}{M})^2 Q(\text{Bou}_M \times \text{Bou}_M)$. Define

$$t = \frac{q}{N} = \frac{NQ(\text{Bou}_M \times \text{Bou}_M)}{M^2}$$

and suppose first that $t$ is an integer. Define a probability measure $\bar{\nu}$ on paths in the skeleton graph by:

- picking an initial state in $\text{Bou}_{N,M}^+$ from distribution $\pi$,
- running the $Q^1$-chain for $t$ steps,
- deleting arrow-labels.

Consider a typical element $(b_0, \to)$ of $\text{Bou}_{N,M}^+$ with $b_0 = (x_0, y_0) + b_0^*$ for some $b_0^* \in \text{Bou}_M$. Recalling definition (12) and using stationarity in the first line below,

$$\begin{aligned}
\bar{\nu}_{(S)}(b_0, b_1) &= t\bar{\nu}\{\text{ paths with first step } b_0 \to b_1\} \\
&= t\pi(b_0, \to)Q^1((b_0, \to), (b_1, \sharp)) \\
&= t\frac{Q_{\text{ent}}(b_0^*)}{q}\frac{Q(b_0^*, b_1^*)}{Q_{\text{ent}}(b_0^*)} \\
&= \frac{t}{q}Q(b_0^*, b_1^*)
\end{aligned}$$

and the same identity holds for general $(b_0, b_1)$. Now define the normalized measure $\nu(\cdot) := N\bar{\nu}(\cdot)$. Then

$$\nu_{(S)}(b_0, b_1) = \frac{tN}{q}Q(b_0^*, b_1^*) = Q(b_0^*, b_1^*)$$

by definition of $t$. Define $\theta_{N,M} = \text{tra}(\nu)$. Properties 3 follow from the constancy properties of the stationary distribution $\pi$ at (13), and assertion (i) of Proposition 4 is immediate from the construction.

Now note that if $t$ is not an integer, then we can apply this construction to $\lfloor t \rfloor$ and to $\lceil t \rceil$, and mix over these two cases.

It remains to prove assertion (ii). Note that the stationary $\pi$ at (13), projected down to $\text{Bou}_{N,M}$, is asymptotically uniform as $N \to \infty$, in the sense that its push-forward under the map $z \to N^{-1}z$ converges weakly to the uniform distribution on the continuous torus $\mathcal{T}^2$.

Note that $\bar{Q}$ defines a Markov transition matrix $Q^*$ on $\text{Bou}_M$ via

$$Q^*(b_0, b_1) = \bar{Q}(b_0, b_1^{\text{reflect}}),$$

where $b_1 \to b_1^{\text{reflect}}$ is the "reflect top boundary with bottom boundary, reflect left boundary with right boundary" procedure above (6).

Define a projection map $\text{pro}_{N,M} : \text{Bou}_{N,M}^+ \to \text{Bou}_M$ as "take modulo $(M, M)$, choosing the boundary point so that the arrow points into the standard $M \times M$ square." So, for instance,

$$\text{pro}_{N,M}((3M + i, 2M - \tfrac{1}{2}), \uparrow) = (i, -\tfrac{1}{2})$$

$$\text{pro}_{N,M}((3M + i, 2M - \tfrac{1}{2}), \downarrow) = (i, M - \tfrac{1}{2}).$$

Here is the key idea in the construction, illustrated in Figure 2. Let $(Y_s, s = 0, 1, 2, \ldots)$ be the $Q^1$-chain, or, more precisely, the chain defined in the same way over the lattice $\mathbb{Z}^2$, so that its values modulo $(N, N)$ are the $Q^1$-chain. Then it is straightforward to verify

$$(X_s := \text{pro}_{N,M}(Y_s), \ s = 0, 1, 2, \ldots) \text{ is the } Q^*\text{-chain}.$$

So the displacement $Y_T - Y_0$ of the $Q^1$-chain over $T$ steps can be written as

$$Y_T - Y_0 = \sum_{s=1}^{T} g(X_{s-1}, X_s); \qquad g(b_0, b_1) = b_1^{\text{reflect}} - b_0.$$

Note the right-hand side does not involve $N$. Suppose first

(14) \qquad the $Q^*$-chain is irreducible.

Then by the strong law of large numbers for additive functionals of a finite Markov chain,

$$T^{-1} \sum_{s=1}^{T} g(X_{s-1}, X_s) \to \bar{g} := \sum_{b,b' \in \text{Bou}_M} \pi^*(b) Q^*(b, b') g(b, b'),$$

where $\pi^*$ is the stationary distribution of $Q^*$. In terms of $t = \frac{NQ(\text{Bou}_M \times \text{Bou}_M)}{M^2}$ steps of the $Q^1$-chain $(Y_s)$, this says that as $N \to \infty$,

$$N^{-1}(Y_t - Y_0) \to \frac{\bar{g} Q(\text{Bou}_M \times \text{Bou}_M)}{M^2} \text{ in probability.}$$

To calculate $\bar{g}$, note that $\pi^*(\cdot)$ is proportional to $Q_{\text{ent}}(\cdot)$. So

$$\bar{g} = \sum_{b,b'} \frac{Q_{\text{ent}}(b)}{Q(\text{Bou}_M \times \text{Bou}_M)} \bar{Q}(b, b')(b' - b)$$

$$= \frac{\sum_{b,b'} (b' - b) \ Q(b, b')}{Q(\text{Bou}_M \times \text{Bou}_M)}$$

$$= \frac{M^2 \, \text{drift}(Q)}{Q(\text{Bou}_M \times \text{Bou}_M)}$$

and so

(15) \qquad $N^{-1}(Y_t - Y_0) \to \text{drift}(Q)$ in probability.

Looking at assertion (ii) of Proposition 4, we noted earlier that the first marginal of $N^{-1}\widetilde{\theta}_{N,M}$ converges weakly to $\mathrm{dist}(U)$, and (15) now shows that the conditional distribution converges to the unit mass at $\mathrm{drift}_1(Q)$, completing the proof. □

Recall the definition (10) of $c_{M,B}$. By compactness, the *inf* is attained, so we have
$$c_{M,B} = \mathbb{E}\,\mathrm{cost}_{M,B}(Q_0, \mathbf{c}),$$
where
$$Q_0 = \int_{\mathcal{Q}_M} Q\,\psi(dQ) \tag{16}$$
for a certain probability measure $\psi$ on $\mathcal{Q}_M$.

PROPOSITION 5. *Given an environment* $\mathbf{c}$ *on the* $N \times N$ *torus, we can define path-flows* $\mu_{N,M}(\cdot|\mathbf{c})$ *with associated flow-volumes* $\mathbf{f}(\cdot|\mathbf{c})$ *on the torus such that:*

(i) $\max_e f(e|\mathbf{c}) \leq B$.
(ii) $\mathbb{E}\,\mathrm{cost}_{(N)}(\mathbf{f}(\cdot|\mathbf{c}), \mathbf{c}) = (\frac{N}{M})^2 c_{M,B}$.
(iii) $\mathrm{tra}(\mu_{N,M}(\cdot|\mathbf{c}))$ *is a measure* $\rho_{N,M}$ *on* $\mathrm{Bou}_{N,M} \times \mathrm{Bou}_{N,M}$ *which does not depend on* $\mathbf{c}$, *and which satisfies Properties 3.*
(iv) *Let* $\widetilde{\rho}_{N,M}$ *be the push-forward of* $\rho_{N,M}$ *under the map* $(z^1, z^2) \to (\frac{1}{N}z^1, \frac{1}{N}z^2)$. *Then as* $N \to \infty$ *with* $M$ *fixed*, $N^{-1}\widetilde{\rho}_{N,M}$ *converges weakly to the uniform probability distribution on* $\mathcal{T}^2 \times \mathcal{T}^2$.

PROOF. Deferring the issue of irreducibility in its hypothesis, Proposition 4 associates with each $Q \in \mathcal{Q}_M$ a measure $\nu^Q$ on skeleton-paths. Define the mixture
$$\nu := \int_{\mathcal{Q}_M} \nu^Q\,\psi(dQ)$$
corresponding to (16). Then $\rho_{N,M} := \mathrm{tra}(\nu)$ satisfies Properties 3 by Proposition 4 and the fact that these properties are preserved under mixtures. Result (iv) follows from Proposition 4(ii) and the definition of isotropic [under $\psi$, the distribution of $\mathrm{drift}_1(Q)$ is uniform on the continuous torus $\mathcal{T}^2$]. The result that
$$\nu_{(S)} = Q_0 \tag{17}$$
follows from Proposition 4(i) and the mixture construction.

Now consider the extended $M \times M$ square. By definition of $c_{M,B}$, there exist path-flows $\mu^0(\cdot|\mathbf{c})$ across the square (depending on the environment $\mathbf{c}$), with flow-volumes $\mathbf{f}^0(\cdot|\mathbf{c})$, such that

$$\max_e f^0(e|\mathbf{c}) \leq B \qquad \forall \mathbf{c}, \tag{18}$$

$$\operatorname{tra}(\mu^0(\cdot|\mathbf{c})) = Q_0 \qquad \forall \mathbf{c}, \tag{19}$$

$$\mathbb{E}\operatorname{cost}_M(\mathbf{f}^0(\cdot|\mathbf{c}), \mathbf{c}) = c_{M,B}. \tag{20}$$

For $b, b' \in \operatorname{Bou}_M$, write $\mu^0_{b,b'}(\cdot|\mathbf{c})$ for the flow $\mu^0(\cdot|\mathbf{c})$ restricted to paths with entrance $b$ and exit $b'$, and normalized to have mass 1.

Recalling the construction of $\nu^Q$ in Proposition 4, we see that $\nu$ can be constructed in 5 steps:

- pick $Q$ from distribution $\psi$,
- pick an initial state $x_0$ in $\operatorname{Bou}^+_{N,M}$ from the stationary distribution of $Q$,
- run the $Q^1$-chain for $t$ steps,
- delete arrow-labels to get a path $(x_0, x_1, \ldots, x_t)$ in the skeleton graph,
- take the probability distribution of the skeleton-paths thus defined, and multiply the measure by $N$.

This gives measure $\nu$. But we can augment the construction in the natural way. Let $\mathbf{c}$ be an environment on the $N \times N$ torus. In each step of the skeleton walk, the successive points $(x_{i-1}, x_i)$ are in the boundary of some subsquare $S_i$. Write $\mathbf{c}^{(i)}$ for the restriction of $\mathbf{c}$ to $S_i$. Expand each step $(x_{i-1}, x_i)$ into a path $(x_{i-1} = y_0, y_1, \ldots, y_m = x_i)$ chosen from distribution $\mu^0_{x_{i-1}, x_i}(\cdot|\mathbf{c}^{(i)})$ (up to translation). Augmenting the construction in this manner gives a measure $\mu_{N,M}(\cdot|\mathbf{c})$ on paths in the $N \times N$ torus. It is easy to check [the key point being that the same transportation measure $Q_0$ appears in (17) and in (19)] that the restriction of these path-flows to a subsquare $S_i$ gives path-flows on $S_i$ agreeing (up to translation) with the path-flows $\mu^0(\cdot|\mathbf{c})$ across the standard $M \times M$ square. So assertions (i) and (ii) follow from (18) and (20).

Finally we outline the technical issue of handling the irreducibility condition in Proposition 4. In the mixture representation (16) let $\Upsilon$ be the push-forward of $\psi$ under the map $Q \to \operatorname{drift}(Q)$. It is not hard to show [cf. Lemma 20(vi)] that $\int_{\mathbb{R}^2} \|u\|_1 \Upsilon(du) < \infty$. Make some particular choice of irreducible transportation measures $\mathbf{Q}^{(u)} \in \mathcal{Q}_M$ with $\operatorname{drift}(\mathbf{Q}^{(u)}) = u$ and construct the mixture

$$Q_0^* = \int_{\mathcal{Q}_M} \mathbf{Q}^{(\operatorname{drift}(Q))} \psi(dQ).$$

Then $\varepsilon Q_0^*$ will be a feasible transportation measure, for sufficiently small $\varepsilon$. So we can repeat the proof above with the mixture (16) replaced by

$$(1-\varepsilon)Q_0 + \varepsilon Q_0^* = \int_{\mathcal{Q}_M} ((1-\varepsilon)Q + \varepsilon \mathbf{Q}^{(\operatorname{drift}(Q))}) \psi(dQ)$$

because now each $(1-\varepsilon)Q + \varepsilon \mathbf{Q}^{(\mathrm{drift}(Q))}$ is irreducible. This gives the same conclusion except that in (ii) $c_{M,B}$ is replaced by $c_{M,B}(\varepsilon) := \mathbb{E}\operatorname{cost}_{M,B}((1-\varepsilon)Q_0 + \varepsilon Q_0^*, \mathbf{c})$. Because $c_{M,B}(\varepsilon) \to c_{M,B}$ as $\varepsilon \downarrow 0$, this conclusion suffices for our later needs. □

3.2. *Patching and smoothing flows.* Here we collect an assortment of simple lemmas, mostly involving flows whose construction does not depend on the environment. Some are used in the next section in proving the upper bound, the others in Section 4.6 in proving the lower bound.

Write $\mathbf{f}^{\mathrm{uni}}$ for the flow-volume on the $N \times N$ torus associated with the standardized global flow which puts equal weight on all minimum-length paths between entrance and exit.

LEMMA 6.  $f^{\mathrm{uni}}(e) = \lfloor N^2/2 \rfloor/(2N^2) \leq \frac{1}{4}$ *for each edge* $e$.

PROOF. Write $\|v-w\|_1$ for the graph-distance (i.e., shortest path length) between vertices $v$ and $w$ on the torus. By symmetry, $f^{\mathrm{uni}}(e)$ does not depend on $e$, so

$$2N^2 f^{\mathrm{uni}}(e) = \sum_e f^{\mathrm{uni}}(e)$$
$$= N^{-3} \sum_v \sum_w \|v-w\|_1$$
$$= N^{-1} \sum_v \|v-(0,0)\|_1$$
$$= 2 \sum_{i=1}^{N-1} \min(i, N-i)$$
$$= \lfloor N^2/2 \rfloor. \qquad \square$$

The next lemma has a similar proof, which we omit.

LEMMA 7. *Let* $U$ *be a uniform random vertex in the* $N \times N$ *torus, and let* $U_{(L)}$ *be a uniform random point in the* $L \times L$ *square centered at* $U$, *where* $1 \leq L \leq N$ *is odd. Then there exists a path-flow* $\mu$ *such that*

$$\operatorname{tra}(\mu) \text{ is the distribution of } (U, U_{(L)})$$

$$\max_e \operatorname{flo}(\mu)(e) \leq \frac{L}{4N^2}.$$

LEMMA 8. *Let $\rho$ be a measure on the boundary $\mathrm{Bou}_M$ of the extended $M \times M$ square, and let $\bar{\rho}$ be the uniform probability distribution on the internal vertices $\mathrm{Bou}_M^o$. There exists a measure $\mu$ on path-space $\Sigma_M$ such that*

$$\mathrm{tra}(\mu) = \rho \times \bar{\rho},$$

$$\max_e \mathrm{flo}(\mu)(e) \leq 2 \max_b \rho(b).$$

PROOF. Use the path from a boundary vertex $b$ to an internal vertex $a$ which makes at most one turn. Check that, for any $e$, the number of pairs $(b,a)$ for which the path traverses $e$ is at most $2M^2$. So

$$f(e) \leq 2M^2 \cdot \max_b \rho(b) \cdot \frac{1}{M^2}. \qquad \square$$

LEMMA 9. *Let $Q$ be a measure on $\mathrm{Bou}_M \times \mathrm{Bou}_M$. There exists a measure $\mu$ on path-space $\Sigma_M$ such that*

$$\mathrm{tra}(\mu) = Q,$$

$$\max_e \mathrm{flo}(\mu)(e) \leq 2\left(\max_b Q_{\mathrm{exi}}(b) + \max_b Q_{\mathrm{ent}}(b)\right).$$

PROOF. For each pair $(b,b') \in \mathrm{Bou}_M \times \mathrm{Bou}_M$, route flow via a uniform internal vertex, and use Lemma 8. $\square$

LEMMA 10. *Let $\theta$ be a measure on the bottom boundary points $\mathrm{Bou}_K^{\mathrm{bottom}}$ of the extended $K \times K$ square, and let $\bar{\theta}$ be the uniform probability measure on the top boundary points $\mathrm{Bou}_K^{\mathrm{top}}$. Then there exists a path-flow $\mu$ such that:*

 (i) $\mathrm{tra}(\mu) = \theta \times \bar{\theta}$;
 (ii) $\max_e \mathrm{flo}(\mu)(e) \leq \max_b \theta(b)$.

REMARK. We use $K$ instead of $M$ here for consistency with its application later.

PROOF OF LEMMA 10. By scaling, we may suppose $\max_b \theta(b) = 1$. Then it is enough to consider the special case $\theta(b) = 1 \,\forall\, b \in \mathrm{Bou}_K^{\mathrm{bottom}}$ and delete surplus flow. In this special case, we route flow from $(x, -\frac{1}{2}) \in \mathrm{Bou}_K^{\mathrm{bottom}}$ to $(y, K - \frac{1}{2}) \in \mathrm{Bou}_K^{\mathrm{top}}$ as follows. If $x = y$, we route straight upward. Otherwise set $d = |y - x| - 1$ and route from $(x, -\frac{1}{2})$, up to $(x, d)$, across to $(y, d)$, up to $(y, K - \frac{1}{2})$.

It is clear that the flow on each vertical edge equals 1. For the horizontal edge $(i-1,d) - (i,d)$, the flow equals

$$\frac{2}{K}|\{(x,y): 0 \leq x \leq i-1, \ i \leq y \leq K-1, y = x+d+1\}|$$

$$\leq \frac{2}{K} \max(i, K-i) \leq 1. \qquad \square$$

The next three lemmas compare the costs of different flows. Note that the same device [the factorization in (21)] is used each time.

LEMMA 11. *For $\frac{1}{4} < B_1 < B_2$,*

$$\mathbb{E} \operatorname{cost}_{(N)}(\mathbf{c}, B_2) \leq \mathbb{E} \operatorname{cost}_{(N)}(\mathbf{c}, B_1)$$

$$\leq \mathbb{E} \operatorname{cost}_{(N)}(\mathbf{c}, B_2) + c^* N^2 \frac{B_2(B_2 - B_1)}{B_2 - 1/4}.$$

PROOF. The left inequality is immediate. Fix $\mathbf{c}$ with $\max_e c(e) \leq c^*$. Let $\mathbf{f}^{(2)}$ be a flow volume for which $\max_e f^{(2)}(e) \leq B_2$ and $\operatorname{cost}(\mathbf{f}^{(2)}, \mathbf{c}) = \operatorname{cost}_{(N)}(\mathbf{c}, B_2)$ and let $\mathbf{f}^{\mathrm{uni}}$ be the uniform flow from Lemma 6. Define $\mathbf{f}^{(1)} = \lambda \mathbf{f}^{\mathrm{uni}} + (1-\lambda) \mathbf{f}^{(2)}$, where $\lambda = \frac{B_2 - B_1}{B_2 - 1/4}$ is the solution of

$$B_1 = \tfrac{1}{4}\lambda + B_2(1-\lambda).$$

Then $\mathbf{f}^{(1)}$ is a standardized global flow with $\max_e f^{(1)}(e) \leq B_1$, and

(21)
$$\begin{aligned}
&\operatorname{cost}(\mathbf{f}^{(1)}, \mathbf{c}) - \operatorname{cost}(\mathbf{f}^{(2)}, \mathbf{c}) \\
&= \sum_e c(e)(f^{(1)}(e) - f^{(2)}(e))(f^{(1)}(e) + f^{(2)}(e)) \\
&\leq c^*(2B_2) \sum_e (\lambda \tfrac{1}{4}) \\
&= c^* B_2 N^2 \lambda,
\end{aligned}$$

establishing the lemma. $\square$

The next lemma relates almost-isotropic transportation measures to isotropic ones. Let $\mu^{\rightarrow}$ be the path-flow across the extended $M \times M$ square which assigns weight 1 to each left-to-right horizontal path; let $Q^{\rightarrow}$ be its transportation measure. Define similarly $\mu^{\leftarrow}$, $\mu^{\uparrow}$, $\mu^{\downarrow}$.

LEMMA 12. *Suppose $Q \in \mathcal{Q}_M$ and*

$$Q + \varepsilon_1 Q^{\rightarrow} + \varepsilon_2 Q^{\leftarrow} + \varepsilon_3 Q^{\uparrow} + \varepsilon_4 Q^{\downarrow} \in \mathcal{Q}_M^{\mathrm{iso}}.$$

*Write $\varepsilon = \varepsilon_1 + \varepsilon_2 + \varepsilon_3 + \varepsilon_4$. Then for $B' \geq B + \varepsilon$,*

$$\mathbb{E}\,\mathrm{cost}_{M,B}(Q,\mathbf{c}) \geq c_{M,B'} - 4c^* B' \varepsilon M^2.$$

PROOF. Given $\mathbf{c}$, let $\mu(\cdot|\mathbf{c})$ be a path-flow attaining the *inf* in the definition (9) of $\mathrm{cost}_{M,B}(Q,\mathbf{c})$. Consider

$$\mu'(\cdot|\mathbf{c}) := \mu(\cdot|\mathbf{c}) + \varepsilon_1 \mu^{\rightarrow} + \varepsilon_2 \mu^{\leftarrow} + \varepsilon_3 \mu^{\uparrow} + \varepsilon_4 \mu^{\downarrow}$$

whose transportation measure is

$$Q' := Q + \varepsilon_1 Q^{\rightarrow} + \varepsilon_2 Q^{\leftarrow} + \varepsilon_3 Q^{\uparrow} + \varepsilon_4 Q^{\downarrow} \in \mathcal{Q}_M^{\mathrm{iso}}.$$

Write $\mathbf{f}$ and $\mathbf{f}'$ for $\mathrm{flo}(\mu(\cdot|\mathbf{c}))$ and $\mathrm{flo}(\mu'(\cdot|\mathbf{c}))$. Note $\max_e f(e) \leq B + \varepsilon \leq B'$ and

$$\mathrm{cost}_M(\mathbf{f}',\mathbf{c}) - \mathrm{cost}_M(\mathbf{f},\mathbf{c}) = \sum_e c(e)(f'(e) + f(e))(f'(e) - f(e))$$

$$\leq \max_e c(e) \times |\mathcal{E}_M|(2B')\varepsilon.$$

So

$$\begin{aligned}
c_{M,B'} &\leq \mathbb{E}\,\mathrm{cost}_{M,B}(Q',\mathbf{c}) \\
&\leq \mathbb{E}\,\mathrm{cost}_M(\mathrm{flo}(\mu'(\cdot|\mathbf{c})),\mathbf{c}) \\
&\leq \mathbb{E}\,\mathrm{cost}_M(\mathrm{flo}(\mu(\cdot|\mathbf{c})),\mathbf{c}) + 2c^* B' |\mathcal{E}_M|\varepsilon \\
&= \mathbb{E}\,\mathrm{cost}_{M,B}(Q,\mathbf{c}) + 4c^* B' M^2 \varepsilon,
\end{aligned}$$

where we have used $|\mathcal{E}_M| = 2M^2$. □

The next lemma relates the cost of an almost-standardized global flow to the cost of a standardized global flow.

LEMMA 13. *Let $\mathbf{f} = \mathbf{f}_1 + \mathbf{f}_2$ be the flow-volume in a standardized global flow on the $N \times N$ torus. Let $0 < \delta < B_0 - \frac{1}{4}$. Suppose*

$$\max_e f_1(e) \leq B_0,$$

$$\max_e f_2(e) \leq \delta.$$

*Then there exists a standardized global flow whose flow-volume $\tilde{\mathbf{f}}$ satisfies*

$$\max_e \tilde{f}(e) \leq B_0 - \delta,$$

$$\mathrm{cost}(\tilde{\mathbf{f}},\mathbf{c}) \leq \mathrm{cost}(\mathbf{f}_1,\mathbf{c}) + \frac{8\delta N^2 B_0^2}{B_0 - 1/4} \max_e c(e)$$

*for every environment $\mathbf{c}$.*

PROOF. Let $\mathbf{f}^{\mathrm{uni}}$ be the "minimum-length" flow-volume from Lemma 6, so $\max_e f^{\mathrm{uni}}(e) \leq 1/4$. For $0 < \lambda < 1$, consider the mixture $\mathbf{f}_\lambda = (1-\lambda)(\mathbf{f}_1 + \mathbf{f}_2) + \lambda \mathbf{f}^{\mathrm{uni}}$ arising from the corresponding mixture of path-flows. Choosing $\lambda = \frac{2\delta}{\delta + B_0 - 1/4}$ which is the solution of

$$B_0 - \delta = (1-\lambda)(B_0 + \delta) + \tfrac{1}{4}\lambda$$

ensures that $\max_e f_\lambda(e) \leq B_0 - \delta$. For any $e$,

$$\begin{aligned}
f_\lambda^2(e) - f_1^2(e) &= (f_\lambda(e) + f_1(e))(f_\lambda(e) - f_1(e)) \\
&\leq \max(0, 2B_0(f_\lambda(e) - f_1(e))) \\
&\leq 2B_0(\delta + \lambda f^{\mathrm{uni}}(e)) \\
&\leq 2B_0\left(\delta + \frac{2\delta}{4B_0 - 1}\right) \\
&\leq \frac{4\delta B_0^2}{B_0 - 1/4}.
\end{aligned}$$

Summing over the $2N^2$ edges $e$ gives the stated bound for $\mathrm{cost}(\tilde{\mathbf{f}}, \mathbf{c})$. □

LEMMA 14. *Identify the edges of the $N \times N$ torus with a subset of the edges of the $(N+1) \times (N+1)$ torus in the natural way, identifying edge $(0,i) - (N-1,i)$ with edge $(0,i) - (N,i)$. Let $\mathbf{c}'$ be an environment on the $(N+1) \times (N+1)$ torus with $\max_e c'(e) \leq c^*$, and let $\mathbf{c}$ be the environment on the $N \times N$ torus induced by the identification. Let $\mathbf{f}$ be a standardized global flow on the $N \times N$ torus such that $\max_e f(e) \leq B$. Then there exists a standardized global flow $\mathbf{f}'$ on the $(N+1) \times (N+1)$ torus such that $\max_e f'(e) \leq B$ and*

$$\mathrm{cost}(\mathbf{f}', \mathbf{c}') \leq \mathrm{cost}(\mathbf{f}, \mathbf{c}) + (N+1)G(c^*, B),$$

*where $G(c^*, B)$ depends only on $c^*$ and $B$.*

PROOF. We construct a flow in three pieces. A path in the $N \times N$ torus induces a path in the $(N+1) \times (N+1)$ torus in a natural way, where necessarily identifying an edge such as $(0,i) - (N-1,i)$ in the former with the two edges $(0,i) - (N,i) - (N-1,i)$ in the latter. So the flow $\mathbf{f}$ induces a flow $\mathbf{f}_1$ with associated path-flow $\mu_1$ on the $(N+1) \times (N+1)$ torus such that, after renormalizing to ensure

$$\mathrm{tra}(\mu_1)(v_1, v_2) = \frac{1}{(N+1)^3}, \quad v_1, v_2 \in [0, N-1]^2,$$

we have

$$\max_e f_1(e) \leq B\left(\frac{N}{N+1}\right)^3, \tag{22}$$

$$\sum_{e \in \mathcal{E}_N} c(e) f_1^2(e) = \left(\frac{N}{N+1}\right)^6 \text{cost}(\mathbf{f}, \mathbf{c}), \tag{23}$$

where $\mathcal{E}_N$ denotes the subset of the edges of the $(N+1) \times (N+1)$ torus identified with the edges of the $N \times N$ torus.

For the second piece, consider the vertex set

$$A := \{(i, N) : 0 \leq i \leq N-1\} \cup \{(N, i) : 0 \leq i \leq N-1\}.$$

We shall define a path-flow $\mu_2$ such that

$$\text{tra}(\mu_2)(v_1, v_2) = (N+1)^{-3}, \qquad (v_1, v_2) \in A \times [0, N-1]^2. \tag{24}$$

We can do this using Lemma 8 with $M = N$, because the boundary points of the extended $N \times N$ square can be mapped via a two-to-one map to $A$. Lemma 8 with $\rho$ uniform, under this mapping, gives a path-flow $\hat{\mu}_2$ such that

$$\text{tra}(\hat{\mu}_2)(v_1, v_2) = 2 \times \frac{1}{4N} \times \frac{1}{N^2}, \qquad (v_1, v_2) \in A \times [0, N-1]^2,$$

$$\max_e \text{flo}(\hat{\mu}_2)(e) \leq 4 \times \frac{1}{4N}. \tag{25}$$

Renormalization gives $\mu_2$ satisfying (24) and

$$\max_e \text{flo}(\mu_2)(e) \leq \frac{2N^2}{(N+1)^3}.$$

For the third piece we want a path-flow $\mu_3$ such that

$$\text{tra}(\mu_3)(v_1, v_2) = \frac{1}{(N+1)^3},$$

$$(v_1, v_2) \in A \times A \cup A \times \{(N, N)\} \cup \{(N, N)\} \times A.$$

We can construct such a flow using only the edges of $\mathcal{E}_{(N+1)} \setminus \mathcal{E}_{(N)}$ which are parallel to the boundary of the $N \times N$ square, that is, of the form $(N, i) - (N, i+1)$ and $(i, N) - (i+1, N)$. It is easy to construct such a flow with

$$\max_e \text{flo}(\mu_3)(e) \leq \frac{1}{N+1}.$$

These edges were not used by the previous flows.

Now superimpose the path-flow $\mu_2$ with its reversal (whose transportation measure is supported on $[0, N-1]^2 \times A$) and with $\mu_3$, and write $\mathbf{f}_2$ for the resulting flow-volume. So we have constructed flow-volumes $\mathbf{f}_1$ and $\mathbf{f}_2$ on the $(N+1) \times (N+1)$ torus such that

$$\mathbf{f}_1 + \mathbf{f}_2 \text{ is a standardized global flow}$$

$$\max_e f_1(e) \leq B \qquad \text{by (22)},$$
$$\max_e f_2(e) \leq \frac{5}{N+1} \qquad \text{by (25)},$$
$$\text{cost}(\mathbf{f}_1, \mathbf{c}') \leq \text{cost}(\mathbf{f}, \mathbf{c}) + 2Nc^*B^2,$$

the latter by (23) and the crude bound for the doubled edges. Applying Lemma 13 to $N+1$ and with $\delta = \frac{5}{N+1}$ establishes the result. □

Finally we state a result, somewhat analogous to the previous lemma, relating $c_{M+1,B}$ to $c_{M,B}$. Because this will be used (see end of Section 4.6) more for aesthetic than essential reasons, we take the opportunity to omit the proof.

LEMMA 15.  $c_{M+1,B} \leq c_{M,B} + (M+1)G'(c^*, B)$, where $G'(c^*, B)$ depends only on $c^*$ and $B$.

3.3. *The upper bound.* Proposition 17 below gives the upper bound in Theorem 2. If the flows in Proposition 5 were exactly standardized global flows, then Proposition 17 would follow immediately. Instead, we will use the following easy smoothing lemma, together with the patching lemmas of Section 3.2.

LEMMA 16.  Let $Z_N$ take values in the $N \times N$ torus. Let $U_L$ be uniform on the $L \times L$ square centered at the origin ($L$ odd), independent of $Z_N$. Suppose $N^{-1}Z_N$ converges in distribution to the uniform probability distribution on $[0,1]^2$. Then, provided $L_N/N \to 0$ sufficiently slowly,

$$\min_{v \in \mathcal{T}_N^2} P(Z_N + U_{L_N} = v)/N^{-2} \to 1.$$

PROPOSITION 17.  There exist standardized global flows $\mathbf{f}(\cdot|\mathbf{c})$ on the $N \times N$ torus such that $\max_e f(e|\mathbf{c}) \leq B$ and such that, for each fixed $M$,

$$\limsup_N N^{-2} \mathbb{E} \, \text{cost}_N(\mathbf{f}(\cdot|\mathbf{c}), \mathbf{c}) \leq M^{-2} c_{M,B}.$$

PROOF. Fix $M|N$ and an environment $\mathbf{c}$ on the $N \times N$ torus. We construct a flow in 4 steps, where only step 2 depends on $\mathbf{c}$. Consider the transportation measure $\rho_{N,M}$ from Proposition 5. The restriction of its entrance marginal to the boundary $\mathrm{Bou}_{(S)}$ of an extended $M \times M$ square $S$ in the natural partition is (up to translation) a measure $\theta$ on $\mathrm{Bou}_M$, which (by Properties 3) does not depend on $S$ and is the same for the exit marginal.

*Step* 1. Each vertex $v$ of the $N \times N$ torus is in some $M \times M$ square $S$ in the natural partition. Construct a path-flow (with some flow-volume $\mathbf{f}_1$) from each $v$ to a random position ($b_0$, say) on the boundary $\mathrm{Bou}_{(S)}$ with distribution proportional to $\theta$, with total volume $N^{-1}$ starting from each $v$.

*Step* 2. Use the path-flow $\mu_{N,M}(\cdot|\mathbf{c})$ in Proposition 5 [with flow-volume $\mathbf{f}_2 = \mathbf{f}_2(\cdot|\mathbf{c})$] to send flow from a typical position $b_0 \in \mathrm{Bou}_{N,M}$ to another position $b_1 \in \mathrm{Bou}_{N,M}$.

*Step* 3. Reverse step 1, sending flow from $b_1$ to a uniform position ($v'$, say) inside an adjacent $M \times M$ square. Write $\mathbf{f}_3$ for the flow-volume.

*Step* 4. Send flow from $v'$ to a uniform random position in the $L_N \times L_N$ square centered at $v'$. Write $\mathbf{f}_4$ for the flow-volume.

Let $\mathbf{f} = \mathbf{f}_1 + \mathbf{f}_2 + \mathbf{f}_3 + \mathbf{f}_4$ be the flow-volume and $\Xi$ the transportation measure for the concatenated flow. Combining Proposition 5(iv) and Lemma 16, provided $L_N/N \to 0$ slowly, we have

$$\min_{y_0, y_1 \in \mathcal{T}_N^2} \Xi(y_0, y_1) = N^{-3} w_N, \tag{26}$$

where $w_N \to 1$ as $N \to \infty$. In step 1, in each square $S$ we are in the setting of Lemma 8, seeking a path-flow whose transportation measure is of the form $\bar{\rho} \times \rho$, where $\rho$ has total mass $N^{-1}M^2$. By Lemma 8, we can choose such a flow to satisfy

$$\max_e f_1(e) \leq 2 \max_b \rho(b) \leq 2N^{-1}M^2.$$

The same bound holds for $\mathbf{f}_3$ in step 3. In step 4, using Lemma 7 and scaling, we can take

$$\max_e f_4(e) \leq N \cdot \frac{L_N}{4N^2} = \frac{L_N}{4N}.$$

So we have

$$\max_e (f_1(e) + f_3(e) + f_4(e)) \leq \frac{4M^2 + (1/4)L_N}{N},$$

while, by Proposition 5,

$$\max_e f_2(e) \leq B.$$

We want to use Lemma 13. Fix $\delta < B - \frac{1}{4}$. Recall (26): by deleting flow, we may assume

$$\Xi(y_0, y_1) = N^{-3} w_N \qquad \forall (y_0, y_1).$$

So $w_N^{-1}\mathbf{f}$ is a standardized global flow. Apply Lemma 13 to $w_N^{-1}\mathbf{f} = w_N^{-1}\mathbf{f}_2 + w_N^{-1}(\mathbf{f}_1 + \mathbf{f}_3 + \mathbf{f}_4)$ with $B_0 = w_N^{-1}B$. For sufficiently large $N$, we have

$$\frac{4M^2 + (1/4)L_N}{Nw_N} < \delta; \qquad \delta < \frac{B}{w_N} - \frac{1}{4}; \qquad \frac{B}{w_N} - \delta < B$$

and so Lemma 13 gives a standardized global flow $\tilde{\mathbf{f}}(\cdot|\mathbf{c})$ such that $\max_e \tilde{f}(e|\mathbf{c}) \leq B$ and

$$\mathrm{cost}_{(N)}(\tilde{\mathbf{f}}(\cdot|\mathbf{c}), \mathbf{c}) \leq w_N^{-2}\mathrm{cost}_{(N)}(\mathbf{f}_2(\cdot|\mathbf{c}), \mathbf{c}) + \frac{8\delta N^2 B^2 w_N^{-2}}{B/w_N - 1/4}c^*.$$

So using Proposition 5(ii),

$$N^{-2}\mathbb{E}\,\mathrm{cost}_{(N)}(\tilde{\mathbf{f}}(\cdot|\mathbf{c}), \mathbf{c}) \leq \frac{c_{M,B}}{M^2 w_N^2} + \frac{8\delta B^2 w_N^{-2} c^*}{B/w_N - 1/4}.$$

Since $w_N \to 1$ and $\delta > 0$ is arbitrary, we see

$$\limsup_N N^{-2}\mathbb{E}\,\mathrm{cost}_N(\tilde{\mathbf{f}}(\cdot|\mathbf{c}), \mathbf{c}) \leq M^{-2} c_{M,B}$$

as $N$ runs through multiples of $M$. Using Lemma 14 to interpolate other values of $N$, we have established Proposition 17. □

**4. The lower bound.** In outline, we will use a concentration inequality (28) and a convexity property (Lemma 18) to derive a rather technical lower bound (Lemma 19) on the cost of any flow across a $M \times M$ square in terms of a smoothed transportation measure across a slightly larger square. We then consider properties of the induced flow across subsquares arising from a standardized global flow (Section 4.4). These ingredients are then combined to prove the lower bound of Proposition 22.

4.1. *A concentration inequality.* Return to the setting (Section 2.4) of the extended $M \times M$ square. Consider environments $\mathbf{c}$ satisfying

(27) $$0 \leq c(e) \leq c^* \qquad \forall e.$$

Given a transportation measure $Q$ on $\mathrm{Bou}_M \times \mathrm{Bou}_M$, and an environment $\mathbf{c}$, we defined

$$\mathrm{cost}_{M,B}(Q, \mathbf{c}) = \inf_\mu \left\{ \mathrm{cost}_M(\mathrm{flo}(\mu), \mathbf{c}) : \mathrm{tra}(\mu) = Q, \max_e \mathrm{flo}(\mu)(e) \leq B \right\}.$$

By compactness, the *inf* is attained. Now consider the effect of changing $\mathbf{c}$ to $\tilde{\mathbf{c}}$ by changing the cost-factor on only a single edge $e_0$. If $\mu$ attains the *inf* above, then

$$\mathrm{cost}_M(\mathrm{flo}(\mu), \tilde{\mathbf{c}}) \leq \mathrm{cost}_M(\mathrm{flo}(\mu), \mathbf{c}) + c^* B^2$$

and it follows that

$$|\cost_{M,B}(Q,\tilde{\mathbf{c}}) - \cost_{M,B}(Q,\mathbf{c})| \leq c^* B^2.$$

So in our model (3) of i.i.d. random environment we can apply the most basic form of *the method of bounded differences* (Azuma–Hoeffding inequality: see, e.g., [18], Section 1.3 or [5], Section 7.2) to conclude that, for $\lambda > 0$,

$$(28)\quad \mathbb{P}(\cost_{M,B}(Q,\mathbf{c}) \leq -\lambda + \mathbb{E}\cost_{M,B}(Q,\mathbf{c})) \leq \exp\left(-\frac{\lambda^2}{2|\mathcal{E}_M|(c^* B^2)^2}\right).$$

Note that $|\mathcal{E}_M| = 2M^2$ [recall convention below definition (10)]. Note also that we only need a one-sided bound.

*Discussion.* The expectation in (28) is order $M^2$, and the inequality says that fluctuations are only of order $M$. If we use the same approach to the kind of "maximal flow subject to i.i.d. edge-capacities" problem mentioned in the Introduction, then we would study a r.v. measuring maximal flow across a $M \times M$ square. As with (28), this r.v. would have fluctuations of order $M$, but now the expectation is also of order $M$, so we do not have an easy concentration result.

Inequality (28) is also the reason we impose the bound $B$ on edge-flows. This bound causes complications in Section 4.3, where to "smooth" transportation measures on the $M \times M$ square we are forced (to preserve the bound $B$) to extend to a larger square.

4.2. *A convexity lemma.* For fixed $\mathbf{c}$ the function (5) $\mathbf{f} \to \cost_M(\mathbf{f},\mathbf{c})$, is convex, so we have the following:

LEMMA 18. *For fixed $\mathbf{c}$, the map $Q \to \cost_{M,B}(Q,\mathbf{c})$ is convex, as a function of feasible transportation measures $Q$ on $\mathrm{Bou}_M \times \mathrm{Bou}_M$.*

4.3. *Smoothing the transportation measure.* Inequality (28) refers to a fixed transportation measure $Q$. We need to use it when $Q$ depends on $\mathbf{c}$. We do this by defining a certain finite set $[\mathcal{S}_\delta(M+2K,K,B)$ below] of transportation measures and then comparing a general $Q$ to a nearby element of that set.

Fix $K|M$. Consider the extended square $[0, M-1]^2$ centered inside the larger extended square $[-K, M+K-1]^2$. Let $\mu$ be a path-flow across $[0, M-1]^2$ which is feasible in the sense

$$\max_e \mathrm{flo}(\mu)(e) \leq B$$

so that $Q := \mathrm{tra}(\mu)$ is a measure on $\mathrm{Bou}_M \times \mathrm{Bou}_M$. We will first use Lemma 10 to construct an *extension* to a path-flow $\mu_{\mathrm{ext}}$ on $[-K, M+K-1]^2$ whose transportation measure $Q_{\mathrm{ext}} := \mathrm{tra}(\mu_{\mathrm{ext}})$ has a smoothness property.

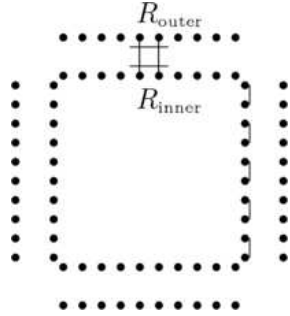

FIG. 3. *Boundary points of an extended $M \times M$ square within matching boundary points of an extended $(M + 2K) \times (M + 2K)$ square. Here $M = 10, K = 2$. A typical $K \times K$ square is shown at the top. The lines on the right side indicate part of the partition of those boundary vertices into blocks of size $K$.*

REMARK. As indicated in Section 2.5, it would be more natural and elegant to do the construction of this section by rerouting flow within the $M \times M$ square. But we are unable to do so. The use of the larger square $[-K, M+K-1]^2$ is purely a technical device, allowing us to create flows with a smooth transportation measure on the larger square without rerouting flow in the inner square.

The construction is illustrated in Figure 3. Partition the $4M$ boundary points into $(4M)/K$ adjacent blocks of size $K$, writing block($b$) for an index of the block containing $b$. Consider a $K \times K$ square $R$ which extends from a block in the boundary of the $M \times M$ square to the parallel block in the boundary of the $(M + 2K) \times (M + 2K)$ square (see Figure 3). Write $R_{\text{inner}}$ and $R_{\text{outer}}$ for the inner and outer boundaries of $R$. Write $\theta_1$ and $\theta_2$ for the restrictions of $Q_{\text{ent}}$ and of $Q_{\text{exi}}$ to $R_{\text{inner}}$. By feasibility, $\theta_1(b) + \theta_2(b) \leq B$ for $b \in R_{\text{inner}}$. Lemma 10 shows that we can create a feasible flow through $R$ whose transportation measure is $\theta_2 \times \bar{\theta}_2 + \bar{\theta}_1 \times \theta_1$, where $\bar{\theta}_i$ is the constant measure on $R_{\text{outer}}$ whose total mass equals $\theta_i(R_{\text{inner}})$. Doing this for each square $R$ constructs a feasible flow $\mu_{\text{smooth}}$ across $[-K, M+K-1]^2$. Because each square $R$ has $2K^2$ edges, and there are $4M/K$ squares, and each edge cost is at most $c^* B^2$, we have

$$\text{cost}_{M+2K}(\text{flo}(\mu_{\text{smooth}}), \mathbf{c}) \leq \text{cost}_M(\text{flo}(\mu), \mathbf{c}) + 8c^* B^2 MK.$$

We emphasize this holds for any path-flow $\mu$ and any environment $\mathbf{c}$ satisfying (27); so the flow may depend on the environment. Note that $Q_{\text{smooth}} := \text{tra}(\mu_{\text{smooth}})$ depends only on $Q := \text{tra}(\mu)$. Recalling the definition (9) of $\text{cost}_{M,B}(Q, \mathbf{c})$, we can minimize over path-flows with prescribed transportation measures to get

(29) $\quad \text{cost}_{M+2K,B}(Q_{\text{smooth}}, \mathbf{c}) \leq \text{cost}_{M,B}(Q, \mathbf{c}) + 8c^* B^2 MK.$

Recall block($b$) denotes the block containing $b$. Define the set $\mathcal{S}(M+2K, K, B)$ to consist of transportation measures $Q$ on $\mathrm{Bou}_{M+2K} \times \mathrm{Bou}_{M+2K}$ with the properties:

(i) $Q(b, b')$ depends only on block($b$) and block($b'$).

(ii) $Q$ is a feasible transportation measure for some path-flow $\mu$ across $[-K, M+K-1]^2$ with $\max_e \mathrm{flo}(\mu)(e) \leq B$.

(iii) $Q(b, b') = 0$ if either $b$ or $b'$ is not one of the $4M$ boundary points which are parallel to the $4M$ boundary points of the $M \times M$ square.

By construction, $Q_{\mathrm{smooth}} := \mathrm{tra}(\mu_{\mathrm{smooth}})$ is an element of $\mathcal{S}(M+2K, K, B)$. We want to relate $\mathrm{drift}(Q_{\mathrm{smooth}})$ to $\mathrm{drift}(Q)$. Write $R$ for a typical $K \times K$ square, $\mathrm{Bou}_{|R}$ for its boundary and $Q_{|R}$ for the transportation measure across $R$, in the construction above. Then

$$(30) \quad \sum_{b,b' \in \mathrm{Bou}_{M+2K}} (b'-b) Q_{\mathrm{smooth}}(b,b') = \sum_{b,b' \in \mathrm{Bou}_M} (b'-b) Q(b,b') + \sum_R \sum_{b,b' \in \mathrm{Bou}_{|R}} (b'-b) Q_{|R}(b,b').$$

For fixed $R$, writing $q^1$ and $q^2$ for the total $Q_{\mathrm{exi}}$-measure and the $Q_{\mathrm{ent}}$-measure of the inner boundary of $R$ (the co-boundary with the $M \times M$ square), we have (using the product form of $Q_{|R}$)

$$\left\| \sum_{b,b' \in \mathrm{Bou}_{|R}} (b'-b) Q_{|R}(b,b') \right\|_1$$
$$\leq |q^1 - q^2| K + (q^1 + q^2)(K/2)$$
$$\leq \tfrac{3}{2} K(q^1 + q^2).$$

Summing over the squares $R$, the $L^1$ size of the second term on the right in (30) is at most $3K\bar{q}$, where $\bar{q} := Q(\mathrm{Bou}_M \times \mathrm{Bou}_M)$. So

$$\|(M+2K)^2 \mathrm{drift}(Q_{\mathrm{smooth}}) - M^2 \mathrm{drift}(Q)\|_1 \leq 3K\bar{q}.$$

Noting that $\|\mathrm{drift}(Q)\|_1 \leq \frac{2\bar{q}}{M}$ and that $1 - \frac{M^2}{(M+2K)^2} \leq 4K/M$, we find

$$\|\mathrm{drift}(Q_{\mathrm{smooth}}) - \mathrm{drift}(Q)\|_1 \leq \frac{4K}{M} \frac{2\bar{q}}{M} + \frac{3K\bar{q}}{(M+2K)^2} \leq \frac{11K\bar{q}}{M^2}.$$

Note that, by feasibility,

$$\bar{q} \leq 2MB.$$

Now let $\delta > 0$ and define $\mathcal{S}_\delta(M+2K, K, B)$ to be the subset of $\mathcal{S}(M+2K, K, B)$ consisting of those $Q$ which satisfy also

(iv) $Q(b, b')$ is a multiple of $\delta/K$.

We first show that the cardinality of $\mathcal{S}_\delta(M+2K,K,B)$ is bounded by a number depending only on $M/K$ and $B/\delta$:

$$(31) \qquad |\mathcal{S}_\delta(M+2K,K,B)| \leq \ell(M/K, B/\delta) := \left(1 + \frac{B}{\delta}\right)^{(4M/K)^2}.$$

To see this, remove constraint (ii) from the definition of $S(M+2K,K,B)$. Then a $Q$ satisfying (i), (iii), (iv) is determined by an array $(m(\beta_1, \beta_2))$ of nonnegative integers, where $\beta_1$ and $\beta_2$ run through the $4M/K$ blocks and

$$Q(b_1, b_2) = m(\beta_1, \beta_2)\delta/K, \qquad b_i \in \beta_i.$$

The latter implies that the volume of flow from $b_1$ to $\beta_2$ is $m(\beta_1, \beta_2)\delta$, and because of the bound $B$ on edge-flows, we must have $m(\beta_1, \beta_2) \leq B/\delta$. This establishes (31).

Define a "trimmed" transportation measure $Q_{\text{trim}} \leq Q_{\text{smooth}}$ by

$Q_{\text{trim}}(b, b')$ is the largest multiple of $\delta$ smaller than $Q_{\text{smooth}}(b, b')$.

So $Q_{\text{trim}}$ is an element of $\mathcal{S}_\delta(M+2K, K, B)$.

Now consider a path-flow $\nu$ with

$$\text{tra}(\nu) = Q_{\text{trim}}; \qquad \max_e \text{flo}(\nu)(e) \leq B.$$

Note that $\widetilde{Q} := Q_{\text{smooth}} - Q_{\text{trim}}$ has $\widetilde{Q}(b_1, b_2) \leq \delta/K$ and so its marginals satisfy $\max(\widetilde{Q}_{\text{ent}}(b), \widetilde{Q}_{\text{exi}}(b)) \leq 4M\delta/K$. Using Lemma 9, we can construct a path-flow $\nu'$ such that

$$\text{tra}(\nu') = \widetilde{Q}; \qquad \max_e \text{flo}(\nu')(e) \leq 16M\delta/K$$

and so

$$\text{tra}(\nu + \nu') = Q_{\text{smooth}}; \qquad \text{flo}(\nu + \nu')(e) \leq \text{flo}(\nu)(e) + 16M\delta/K \leq B',$$

where we set $B' := B + 16M\delta/K$. Because increasing the flow-volume across an edge by $\eta$ (subject to the bound $B'$) can increase the cost associated with that edge by at most $2c^*B'\eta$, it follows that

$$\text{cost}_{M+2K}(\text{flo}(\nu + \nu'), \mathbf{c}) \leq \text{cost}_{M+2K}(\text{flo}(\nu), \mathbf{c}) + 2B'c^*(16M\delta/K)|\mathcal{E}_{M+2K}|.$$

Minimizing over path-flows with prescribed transportation measures, we see

$$(32) \quad \begin{aligned} &\text{cost}_{M+2K, B'}(Q_{\text{smooth}}, \mathbf{c}) \\ &\leq \text{cost}_{M+2K, B}(Q_{\text{trim}}, \mathbf{c}) + 64B'c^*\delta M(M+2K)^2/K. \end{aligned}$$

We now bring into play the probability measure on environments. Write

$$(33) \qquad m_B(Q) := \mathbb{E}\,\text{cost}_{M+2K, B}(Q, \mathbf{c}).$$

Consider, for each $\mathbf{c}$, some feasible path-flow $\mu(\cdot|\mathbf{c})$ across the extended $M \times M$ square (this notation emphasizes we allow dependence on $\mathbf{c}$); write $Q_{\text{smooth}}(\cdot|\mathbf{c})$ and $Q_{\text{trim}}(\cdot|\mathbf{c})$ for its smoothed and trimmed transportation measures. From the concentration inequality (28) applied to $M + 2K$, and the finiteness of $\mathcal{S}_\delta(M + 2K, K, B)$,

$$\mathbb{P}(\text{cost}_{M+2K,B}(Q_{\text{trim}}(\cdot|\mathbf{c}), \mathbf{c}) \leq -\lambda + m_B(Q_{\text{trim}}(\cdot|\mathbf{c})))$$
$$\leq \ell(M/K, B/\delta) \exp\left(-\frac{\lambda^2}{(2(M+2K)c^* B^2)^2}\right).$$

This remains true if we replace the first occurrence of $Q_{\text{trim}}$ by the larger $Q_{\text{smooth}}$. Combining with (29) and rescaling $\lambda$,

$$\mathbb{P}(\text{cost}_{M,B}(Q(\cdot|\mathbf{c}), \mathbf{c}) \leq -(8+\lambda)c^* B^2 MK + m_B(Q_{\text{trim}}(\cdot|\mathbf{c})))$$
$$\leq \ell(M/K, B/\delta) \exp\left(-\frac{\lambda^2 M^2 K^2}{2(M+2K)^2}\right).$$

We tidy slightly by setting $K = M/L$ for some $L \geq 4$ and setting $\lambda = 4$ to get

$$\mathbb{P}(\text{cost}_{M,B}(Q(\cdot|\mathbf{c}), \mathbf{c}) \leq -12c^* B^2 M^2/L + m_B(Q_{\text{trim}}(\cdot|\mathbf{c})))$$
$$\leq \ell(L, B/\delta) \exp(-M^2/L^2).$$

Note that for any feasible $Q$ there is an upper bound on $m_B(Q)$ implied by feasibility, which works out as $m_B(Q_{\text{trim}}(\cdot|\mathbf{c})) \leq B^2 c^* 2(M+2K)^2$. Applying the elementary fact, for nonnegative r.v.'s $X, Y$,

if $P(X \leq Y - a) \leq b$ and $Y \leq y_0$, then $EX \geq EY - a - by_0$

shows

$$\mathbb{E}\,\text{cost}_{M,B}(Q(\cdot|\mathbf{c}), \mathbf{c}) \geq \mathbb{E} m_B(Q_{\text{trim}}(\cdot|\mathbf{c})) - 12c^* B^2 M^2/L$$
$$- B^2 c^* 2(M+2K)^2 \ell(L, B/\delta) \exp(-M^2/L^2).$$

Now (32) says

$$m_B(Q_{\text{trim}}(\cdot|\mathbf{c}), \mathbf{c}) \geq m_{B'}(Q_{\text{smooth}}(\cdot|\mathbf{c}), \mathbf{c}) - 64B' c^* \delta M(M+2K)^2/K.$$

Convexity (Lemma 18) and Jensen's inequality imply

$$\mathbb{E} m_{B'}(Q_{\text{smooth}}(\cdot|\mathbf{c})) \geq m_{B'}(\mathbb{E} Q_{\text{smooth}}(\cdot|\mathbf{c})).$$

Combining the last three inequalities and using the fact $M + 2K \leq \frac{3}{2}M$ (because $L \geq 4$),

$$\mathbb{E}\,\text{cost}_{M,B}(Q(\cdot|\mathbf{c}), \mathbf{c}) \geq m_{B'}(\mathbb{E} Q_{\text{smooth}}(\cdot|\mathbf{c})) - 12 c^* B^2 \frac{M^2}{L}$$
$$- 144 B' c^* \delta M^2 L - 5 c^* B^2 M^2 \ell(L, B/\delta) \exp(-M^2/L^2).$$

We can summarize the results of this section as the next lemma:

LEMMA 19. *Let $M = KL$, where $L \geq 4$. For each environment $\mathbf{c}$ with $\max_e c(e) \leq c^*$, let $\mu(\cdot|\mathbf{c})$ be a path-flow across the extended $M \times M$ square with $\max_e \mathrm{flo}(\mu(\cdot|\mathbf{c}))(e) \leq B$. Let $Q(\cdot|\mathbf{c})$ be its transportation measure and $Q_{\mathrm{smooth}}(\cdot|\mathbf{c})$ the smoothed transportation measure across the extended $(M+2K) \times (M+2K)$ square. Then*

$$\|\mathrm{drift}(Q_{\mathrm{smooth}}(\cdot|\mathbf{c})) - \mathrm{drift}(Q(\cdot|\mathbf{c}))\|_1 \leq \frac{22KB}{M} \tag{34}$$

*and for any $\delta > 0$,*

$$\begin{aligned}
M^{-2} &\mathbb{E}\,\mathrm{cost}_M(\mathrm{flo}(\mu(\cdot|\mathbf{c})), \mathbf{c}) \\
&\geq M^{-2} m_{B'}(\mathbb{E} Q_{\mathrm{smooth}}(\cdot|\mathbf{c})) - 12 c^* B^2/L \\
&\quad - 144 B' c^* \delta L - 5 c^* B^2 \ell(L, B/\delta) \exp(-M^2/L^2),
\end{aligned} \tag{35}$$

*where $B' = B + 16\delta L$.*

**4.4. Local flows induced by a global flow.** This section derives identities (Lemma 20) which restate certain quantities for flows on the $N \times N$ torus in terms of induced empirical distributions across $M \times M$ squares. Of course, this is "just bookkeeping"; the point is to obtain quantities defined on spaces depending on $M$ not on $N$, so that it makes sense later to let $N \to \infty$. Until further notice, fix $2 \leq M < N$, a global environment $\hat{\mathbf{c}}$ on the $N \times N$ torus and a standardized global flow $\mu$. We start with some observations.

OBSERVATION 1. A path $\pi$ in the $N \times N$ torus from $z^1$ to $z^2$ induces (by using the same increments) a path in the lattice $\mathbb{Z}^2$ from $z^1$ to a point $z^*$ of the form $z^* = z^2 + (iN, jN)$, where the integer pair $(i, j)$ indicates net winding around the torus. Define

$$\omega_N(\pi) = (z^* - z^1)/N \in \mathbb{R}^2 \tag{36}$$

so that $\omega_N(\pi)$ indicates both the normalized relative positions of origin and destination of the path within the torus, and also the number of windings.

OBSERVATION 2. Write $\mathcal{M}(S)$ for the space of measures on a space $S$. If $\psi$ is a probability distribution on $\mathcal{M}(S)$, then (under the natural integrability condition) it has a mean measure $\bar{\psi} \in \mathcal{M}(S)$:

$$\bar{\psi}(\cdot) = \int_{\mathcal{M}(S)} \nu(\cdot)\,\psi(d\nu).$$

OBSERVATION 3. Given a measure $\nu \in \mathcal{M}(S \times S')$ with first marginal $\nu_1 \in \mathcal{M}(S)$, there is for $s \in S$ a conditional (probability) distribution $\nu(\cdot|s) \in \mathcal{M}(S')$ satisfying

$$\nu(A \times B) = \int_A \nu(B|s)\nu_1(ds).$$

Recall that $\Sigma_M$ and $\mathcal{C}_M$ denote the sets of paths $\sigma$ and of environments $\mathbf{c}$ on the extended $M \times M$ square. Given $(x,y) \in [0, N-1]^2$, consider $[x, x+M-1] \times [y, y+M-1]$ as a subsquare of the $N \times N$ torus. The global environment $\hat{\mathbf{c}}$ induces a local environment on the subsquare which, by translation, becomes an element $\mathbf{c}_{x,y} = (c_{x,y}(e))$ of $\mathcal{C}_M$:

$$\tag{37} c_{x,y}(e) = \hat{c}(e + (x,y)), \qquad e \in \mathcal{E}_M,$$

where $e + (x,y)$ denotes the edge $e$ translated by $(x,y)$.

Consider a path $\pi$ in the $N \times N$ torus. Then $\pi$ intersects the subsquare $[x, x+M-1] \times [y, y+M-1]$ via some number (maybe zero) of segments $\sigma$, which (by translation) can be regarded as elements of $\Sigma_M$. Call these segments $\sigma_i(\pi; x, y)$, say. So the standardized global flow $\mu$ induces a measure $\mu_{x,y}$ on $\Sigma_M$ via

$$\mu_{x,y}(\cdot) = \sum_i \mu\{\pi : \sigma_i(\pi; x, y) \in \cdot\}.$$

We need to record also the value of $\omega_N(\pi)$ from (36), which leads us to define the measure $\mu_{x,y}^+$ on $\Sigma_M \times \mathbb{R}^2$:

$$\mu_{x,y}^+(\cdot \times \cdot) = \sum_i \mu\{\pi : \sigma_i(\pi; x, y) \in \cdot, \omega_N(\pi) \in \cdot\}.$$

Note

$$(\mathbf{c}_{x,y}, \mu_{x,y}^+) \text{ is an element of } \mathcal{C}_M \times \mathcal{M}(\Sigma_M \times \mathbb{R}^2).$$

Now take $(x,y)$ as a uniform random vertex of the $N \times N$ torus, so that the random element $(\mathbf{c}_{x,y}, \mu_{x,y}^+)$ has some probability distribution, which we call $\Gamma$. So

$$\tag{38} \Gamma \text{ is a probability measure on } \mathcal{C}_M \times \mathcal{M}(\Sigma_M \times \mathbb{R}^2).$$

In words, $\Gamma$ describes the empirical distribution of local (i.e., across $M \times M$ squares) flows, jointly with the local environment and the relative source-destination vector of flow-paths. As indicated initially, $\Gamma$ contains quite a lot of information about the global flow and environment, $\mu$ and $\hat{\mathbf{c}}$, as we next describe.

*The identity for costs.* Write $\nu^+ \to \nu$ for the "take first marginal" map $\mathcal{M}(\Sigma_M \times \mathbb{R}^2) \to \mathcal{M}(\Sigma_M)$. Then

$$(\mathbf{c}, \nu^+) \to \operatorname{cost}_M(\operatorname{flo}(\nu), \mathbf{c})$$

is a functional defined on $\mathcal{C}_M \times \mathcal{M}(\Sigma_M \times \mathbb{R}^2)$. We assert

$$\tag{39} M^{-2} \int_{\mathcal{C}_M \times \mathcal{M}(\Sigma_M \times \mathbb{R}^2)} \operatorname{cost}_M(\operatorname{flo}(\nu), \mathbf{c})\, \Gamma(d\mathbf{c}, d\nu^+) = N^{-2} \operatorname{cost}(\operatorname{flo}(\mu), \hat{\mathbf{c}}).$$

This holds because the integral equals
$$N^{-2} \sum_{(x,y)} \mathrm{cost}_M(\mathrm{flo}(\mu_{x,y}), \mathbf{c}_{x,y})$$
and the sum here equals $M^2 \mathrm{cost}(\mathrm{flo}(\mu), \hat{\mathbf{c}})$ because for each edge $e \in \mathcal{E}_{(N)}$, the cost $\hat{c}(e)f^2(e)$ is counted in exactly $M^2$ subsquares.

*The identity for drift.* As in Observation 1, consider a path $\pi$ in the discrete torus $[0, N-1]^2$ from $z^1$ to $z^2$ and the induced path in the lattice $\mathbb{Z}^2$ from $z^1$ to a point $z^*$. Consider the natural partition $\Pi_0$ of $[0, N-1]^2$ into $N^2/M^2$ subsquares $S_{x,y} = [x, x+M-1] \times [y, y+M-1]$. Then
$$z^* - z^1 = \sum_{S_{x,y} \in \Pi_0} \sum_\sigma (\mathrm{exi}(\sigma) - \mathrm{ent}(\sigma))$$
where the second sum is over the path fragments $\sigma = \sigma_i(\pi; x, y)$, where $\pi$ intersects $S_{x,y}$. Averaging over shifts of the partition gives

(40) $$z^* - z^1 = M^{-2} \sum_{(x,y) \in [0,N-1]^2} \sum_\sigma (\mathrm{exi}(\sigma) - \mathrm{ent}(\sigma)).$$

Now let $\Gamma_2$ be the second marginal of $\Gamma$. So $\Gamma_2$ is a probability distribution on $\mathcal{M}(\Sigma_M \times \mathbb{R}^2)$. Let $\bar{\Gamma}_2$ be its mean measure (Observation 2). Concretely, $\bar{\Gamma}_2 = N^{-2} \sum_{x,y} \mu^+_{x,y}$. For $z \in \mathbb{Z}^2$,
$$\bar{\Gamma}_2\left(\cdot \times \left\{\frac{z}{N}\right\}\right) = N^{-2} \sum_{x,y} \mu^+_{x,y}\left(\cdot \times \left\{\frac{z}{N}\right\}\right).$$

Using linearity of the map $\nu \to \mathrm{drift}(\mathrm{tra}(\nu))$,

(41) 
$$\mathrm{drift}\left(\mathrm{tra}\left(\bar{\Gamma}_2\left(\cdot \times \frac{z}{N}\right)\right)\right)$$
$$= M^{-2} N^{-2} \sum_{x,y} \int \sum_i (\mathrm{exi}(\sigma_i(\pi; x, y))$$
$$- \mathrm{ent}(\sigma_i(\pi; x, y)))1_{(\omega_N(\pi)=z/N)}\mu(d\pi)$$
$$= N^{-2} z \mu\left\{\pi : \omega_N(\pi) = \frac{z}{N}\right\},$$

the last equality by (40). Write
$$\Upsilon_N(\cdot) = N^{-1} \mu\{\pi : \omega_N(\pi) \in \cdot\}.$$

From the definition of *standardized global flow*,

(42) the push-forward of $\Upsilon_N$ under the map $u \to u \bmod (1,1)$ is the uniform probability distribution on $\{N^{-1} z : z \in [0, N-1]^2\}$.

Note that (41) can be rewritten as

(43) $$\mathrm{drift}(\mathrm{tra}(\bar{\Gamma}_2(\cdot \times \{u\}))) = u \Upsilon_N(u), \qquad u \in \mathbb{R}^2.$$

*Identity for marginals of transportation measure.* When a segment $\sigma_i(\pi; x, y)$ exits the square $[x, x+M-1] \times [y, y+M-1]$ via some boundary vertex, it must enter an adjacent square. This fact easily implies the following identity. Given $u = z/N \in \mathbb{R}^2$, write $Q = \mathrm{tra}(\Gamma_2(\cdot \times \{u\}))$. Then

(44) the push-forward of $Q_{\mathrm{exi}}$ (restricted to $\mathrm{Bou}_M$)

under the map $b \to b^{\mathrm{reflect}}$

equals $Q_{\mathrm{ent}}$ (restricted to $\mathrm{Bou}_M$).

*Identity for terminal vertices.* For a path $\sigma \in \Sigma_M$, let $\mathrm{inte}_M(\sigma) \in \{0, 1, 2\}$ be the number of endpoints [$\mathrm{ent}(\sigma)$ and $\mathrm{exi}(\sigma)$] which are in the interior $\mathcal{B}_M^o$ rather than the boundary $\mathrm{Bou}_M$. Here *inte* is a mnemonic for *interior*. Because each path in $\Sigma_{(N)}$ has exactly two endpoints, and an endpoint appears in exactly $M^2$ subsquares,

$$\sum_{x,y} \int_{\Sigma_M} \mathrm{inte}_M(\sigma) \mu_{x,y}(d\sigma) = 2M^2 \mu(\Sigma_{(N)}) = 2M^2 N.$$

Writing $\bar{\Gamma}_{21}(\cdot) = N^{-2} \sum_{x,y} \mu_{x,y}(\cdot)$, so that $\bar{\Gamma}_{21}$ is the marginal of $\bar{\Gamma}_2$ above, we find

(45) $$\int_{\Sigma_M} \mathrm{inte}_M(\sigma) \, \bar{\Gamma}_{21}(d\sigma) = 2M^2/N.$$

*Bound for windings.* For a path $\pi$ in the $N \times N$ torus,

$$\sum_{e \in \mathcal{E}_{(N)}} n(\pi, e) \geq N \|\omega_N(\pi)\|_1$$

and so for a standardized global flow $\mu$ with flow-volume $\mathbf{f}$,

$$\sum_{e \in \mathcal{E}_{(N)}} f(e) \geq N \int \|\omega_N(\pi)\|_1 \, \mu(d\pi).$$

If the capacity constraint $\max_e f(e) \leq B$ is satisfied, then the left-hand side is at most $2N^2 B$ and so

(46) $$\int \|\omega_N(\pi)\|_1 \, \mu(d\pi) \leq 2NB.$$

*Identity for total mass.* From the definitions we can write the total mass of $\bar{\Gamma}_2$ as

$$\bar{\Gamma}_2(\Sigma_M \times \mathbb{R}^2) = N^{-2} \sum_{x,y} \int H_{x,y}(\pi) \mu(d\pi),$$

where $H_{x,y}(\pi)$ is the number of segments of $\pi$ which intersect the square $[x, x + M - 1] \times [y, y + M - 1]$. Because each step enters exactly $M$ such squares, and $\text{ent}(\pi)$ is in exactly $M^2$ such squares,

$$\bar{\Gamma}_2(\Sigma_M \times \mathbb{R}^2) = N^{-2}\Big(M^2 N + M \int \text{len}(\pi)\,\mu(d\pi)\Big), \tag{47}$$

where $\text{len}(\pi)$ is the length of the path $\pi$ and where the factor $N$ is the total mass of $\mu$.

So far we have worked with a fixed standardized global flow $\mu$ and a fixed global environment $\hat{\mathbf{c}}$, so rewrite the $\Gamma$ at (38) as $\Gamma^{\mu,\hat{\mathbf{c}}}$ to indicate this explicitly. Now suppose the environment $\hat{\mathbf{c}}$ is random according to our probability model (3), let $\mu(\cdot|\hat{\mathbf{c}})$ be a standardized global flow depending on the realization of $\hat{\mathbf{c}}$, and write

$$\Gamma(\cdot) = \mathbb{E}\Gamma^{\mu(\cdot|\hat{\mathbf{c}}),\hat{\mathbf{c}}}(\cdot) \tag{48}$$

for the mixed probability distribution on $\mathcal{C}_M \times \mathcal{M}(\Sigma_M \times \mathbb{R}^2)$. Note that the marginal ($\Gamma_1$, say) on $\mathcal{C}_M$ is clearly just the i.i.d. distribution (3). Note also that (39) implies

$$\begin{aligned}M^{-2} \int_{\mathcal{C}_M \times \mathcal{M}(\Sigma_M \times \mathbb{R}^2)} &\text{cost}_M(\text{flo}(\nu), \mathbf{c})\,\Gamma(d\mathbf{c}, d\nu^+) \\ &= N^{-2} \mathbb{E}\,\text{cost}(\text{flo}(\mu(\cdot|\hat{\mathbf{c}})), \hat{\mathbf{c}}).\end{aligned} \tag{49}$$

We now want to replace the probability measure $\Gamma$ by a (nonprobability) measure $\bar{\Gamma}$ on $\mathcal{C}_M \times \Sigma_M \times \mathbb{R}^2$ defined (loosely speaking) by replacing the conditional distribution (given environment) on path-flows by its conditional mean. Precisely, for $\mathbf{c} \in \mathcal{C}_M$, the probability distribution $\Gamma$ has an associated conditional distribution $\Gamma_\mathbf{c}$ on $\mathcal{M}(\Sigma_m \times \mathbb{R}^2)$. Define $\bar{\Gamma}_\mathbf{c} \in \mathcal{M}(\Sigma_m \times \mathbb{R}^2)$ as the mean measure (Observation 2) of $\Gamma_\mathbf{c}$, and then define

$$\bar{\Gamma}(d\mathbf{c}, \cdot, \cdot) = \Gamma_1(d\mathbf{c})\bar{\Gamma}_\mathbf{c}(\cdot, \cdot). \tag{50}$$

Lemma 20 gives identities and inequalities relating $\bar{\Gamma}$ to the underlying global flows in the random environment. Identities (43)–(45) and bound (46), which do not involve interaction between flows and environment, extend to the random environment setting in obvious ways, recorded as (iii)–(vi) below. We handle (49) as follows. The map

$$\nu^+ \to \nu \to \text{flo}(\nu) \quad \text{from } \mathcal{M}(\Sigma_M \times \mathbb{R}^2) \to \mathcal{M}(\Sigma_M) \to \mathcal{F}_M \tag{51}$$

takes $\bar{\Gamma}_\mathbf{c}$ to some element of $\mathcal{F}_M$, say, $\mathbf{f}_\mathbf{c}$. Applying Jensen's inequality to the convex function $\nu^+ \to \text{cost}_M(\text{flo}(\nu), \mathbf{c})$ and the probability measure $\Gamma_\mathbf{c}$ gives

$$\text{cost}_M(\mathbf{f}_\mathbf{c}, \mathbf{c}) \leq \int_{\mathcal{M}(\Sigma_M \times \mathbb{R}^2)} \text{cost}_M(\text{flo}(\nu), \mathbf{c})\,\Gamma_\mathbf{c}(d\nu^+).$$

Integrating both sides against $\Gamma_1(d\mathbf{c})$ and using (49) gives the inequality stated in (ii) below.

LEMMA 20. *For arbitrary standardized global flows $\mu(\cdot|\hat{\mathbf{c}})$ on the $N \times N$ torus, define the measure $\bar{\Gamma}$ on $\mathcal{C}_M \times \Sigma_M \times \mathbb{R}^2$ by (48), (50):*

(i) $\bar{\Gamma}(d\mathbf{c}, \cdot, \cdot) = \Xi_1(d\mathbf{c}) \, \bar{\Gamma}_{\mathbf{c}}(\cdot, \cdot)$, *where $\Xi_1$ is the i.i.d. probability distribution* (3) *on $\mathcal{C}_M$ and where $\bar{\Gamma}_{\mathbf{c}} \in \mathcal{M}(\Sigma_M \times \mathbb{R}^2)$ for each $\mathbf{c} \in \mathcal{C}_m$.*

(ii) *For $\mathbf{f}_{\mathbf{c}}$ defined below* (51),

$$M^{-2}\mathbb{E}\,\mathrm{cost}_M(\mathbf{f}_{\mathbf{c}}, \mathbf{c}) \leq N^{-2}\mathbb{E}\,\mathrm{cost}(\mathrm{flo}(\mu(\cdot|\hat{\mathbf{c}})), \hat{\mathbf{c}}).$$

(iii) *Write $\bar{\Gamma}_2(\cdot, \cdot) = \bar{\Gamma}(\mathcal{C}_M \times \cdot \times \cdot)$ for the marginal measure on $\Sigma_M \times \mathbb{R}^2$. For $D \subset \mathbb{R}^2$, write $Q^{[D]} = \mathrm{tra}(\bar{\Gamma}_2(\cdot \times D))$. Then*

$$\mathrm{drift}(Q^{[D]}) = \int_D u\,\Upsilon_N(du),$$

*where*

$$\Upsilon_N(\cdot) = N^{-1}\mathbb{E}\mu\{\pi : \omega_N(\pi) \in \cdot|\hat{\mathbf{c}}\}.$$

(iv) *For $D \subset \mathbb{R}^2$,*

*the push-forward of $Q^{[D]}_{exi}$ (restricted to $\mathrm{Bou}_M$)*

*under the map $b \to b^{\mathrm{reflect}}$*

*equals $Q^{[D]}_{ent}$ (restricted to $\mathrm{Bou}_M$).*

(v) *For $\bar{\Gamma}_{21}(\cdot) = \bar{\Gamma}_1(\cdot \times \mathbb{R}^2)$,*

$$\int_{\Sigma_M} \mathrm{inte}_M(\sigma)\,\bar{\Gamma}_{21}(d\sigma) = 2M^2/N.$$

*For the remaining parts assume $\max_{e \in \mathcal{E}_{(N)}} \mathrm{flo}(\mu(\cdot|\hat{\mathbf{c}}))(e) \leq B\,\forall \hat{\mathbf{c}}$.*

(vi) $\int_{\mathbb{R}^2} \|u\|_1\,\Upsilon_N(du) \leq 2B$.
(vii) $\max_{e \in \mathcal{E}_M} f_{\mathbf{c}}(e) \leq B\,\forall \mathbf{c}$.
(viii) $\bar{\Gamma}(\mathcal{C}_M \times \Sigma_M \times \mathbb{R}^2) \leq M^2 N^{-1} + 2MB$.

PROOF. Part (i) holds by definition and parts (ii)–(vi) were discussed above. Part (vii) is straightforward. For part (viii), by (47),

$$\bar{\Gamma}(\mathcal{C}_M \times \Sigma_M \times \mathbb{R}^2) = M^2 N^{-1} + N^{-2}M\mathbb{E}\int \mathrm{len}(\pi)\,\mu(d\pi|\hat{\mathbf{c}}).$$

Now the integral equals

$$\sum_{e \in \mathcal{E}_{(N)}} \mathrm{flo}(\mu(\cdot|\hat{\mathbf{c}}))(e) \leq B|\mathcal{E}_{(N)}|,$$

giving (viii). $\square$

4.5. *Properties of weak limits.* Lemma 20 dealt with arbitrary standardized global flows $\mu(\cdot|\hat{\mathbf{c}})$. So we can apply it to flows attaining the *inf* in the definition

$$\mathrm{cost}_{(N)}(\mathbf{c}, B) := \inf\{\mathrm{cost}(\mathbf{f}, \mathbf{c}) : \mathbf{f} \text{ a standardized global flow satisfying (4)}\}$$

and then see what happens in the $N \to \infty$ limit.

Define

(52) $$\gamma_*(\kappa, B) := \liminf_N N^{-2} \mathbb{E}\,\mathrm{cost}_{(N)}(\mathbf{c}, B).$$

Let $\Sigma_M^+ \subset \Sigma_M$ be the subset of paths $\sigma$ which go *across* the $M \times M$ square in the sense that both $\mathrm{ent}(\sigma)$ and $\mathrm{exi}(\sigma)$ are in the boundary $\mathrm{Bou}_M$.

LEMMA 21. *Fix $M \geq 2$. There exist measures $\Xi_{\mathbf{c}}(\cdot, \cdot)$, $\mathbf{c} \in \mathcal{C}_M$ on $\Sigma_M^+ \times \mathbb{R}^2$ with the following properties. Let $\mathbf{f_c}$ be the flow-volume for the path-flow $\Xi_{\mathbf{c}}(\cdot \times \mathbb{R}^2)$:*

(53) $$\max_e f_{\mathbf{c}}(e) \leq B \qquad \forall \mathbf{c},$$

(54) $$M^{-2} \mathbb{E}\,\mathrm{cost}_M(\mathbf{f_c}, \mathbf{c}) \leq \gamma_*(\kappa, B).$$

*For bounded $D \subset \mathbb{R}^2$, write $Q^{[D]} = \mathrm{tra}(\mathbb{E}\Xi_{\mathbf{c}}(\cdot \times D))$. Then*

(55) $$Q^{[D]} \in \mathcal{Q}_M,$$

(56) $$\mathrm{drift}(Q^{[D]}) = \int_D u\,\Upsilon(du),$$

*where $\Upsilon$ is a probability measure on $\mathbb{R}^2$ whose push-forward under the map $u \to u \bmod (1,1)$ is the uniform probability distribution on $\mathcal{T}^2$.*

Note that the final sentence implies (via a Radon–Nikodym argument) that $Q := \mathrm{tra}(\mathbb{E}\Xi_{\mathbf{c}}(\cdot \times \mathbb{R}^2))$ has a "density representation"

$$Q = \int_{\mathbb{R}^2} Q^{[u]}\,\Upsilon(du),$$

where

$$Q^{[u]} \in \mathcal{Q}_M; \qquad \mathrm{drift}(Q^{[u]}) = u,$$

which implies $Q \in \mathcal{Q}_M^{\mathrm{iso}}$. This result is the fundamental motivation for considering $\mathcal{Q}_M^{\mathrm{iso}}$.

PROOF OF LEMMA 21. Recall the vague topology on $\mathcal{M}(\mathcal{C}_M \times \Sigma_M \times \mathbb{R}^2)$, that is, convergence of integrals of continuous functions with compact support. Recall that in this topology the condition for relative compactness of $(\Gamma_N)$ is

$$\sup_N \Gamma_N(K) < \infty \qquad \forall K \text{ compact.}$$

Recall also the Fatou-like result:

$$\text{if } \Gamma_N \to \Xi \text{ vaguely,} \qquad \text{then } \int h \, d\Xi \leq \liminf_N \int h \, d\Gamma_N$$

for continuous $h \geq 0$.

For each $N$, we consider the standardized global flows attaining the *inf* in the definition of $\text{cost}_{(N)}(\mathbf{c}, B)$ and then define $\bar{\Gamma}_N$ to be the measure $\bar{\Gamma}$ associated with these flows as in Lemma 20. The relative compactness condition holds by (viii) of Lemma 20, so we can first take a subsequence of $N$ where the lim inf in (52) obtains, and then a subsequence such that $\bar{\Gamma}_N \to \Xi$ weakly, for some measure $\Xi$ on $\mathcal{C}_M \times \Sigma_M \times \mathbb{R}^2$. By (v) of Lemma 20,

$$\int \text{inte}_M(\sigma) \, \Xi(d\mathbf{c}, d\sigma, du) = 0,$$

which says that $\Xi$ is supported on $\mathcal{C}_M \times \Sigma_M^+ \times \mathbb{R}^2$, that is, we may replace $\Sigma_M$ by $\Sigma_M^+$. For each $N$, the marginal of $\Gamma_N$ on $\mathcal{C}_M$ is the i.i.d. law $\Xi_1$, so the marginal of $\Xi$ is the same law $\Xi_1$, and so there is a representation

$$\Xi(d\mathbf{c}, \cdot, \cdot) = \Xi_1(d\mathbf{c}) \, \Xi_{\mathbf{c}}(\cdot, \cdot)$$

which serves to define $\Xi_{\mathbf{c}}(\cdot, \cdot)$. Properties (53)–(56) now follow as $N \to \infty$ limits of the corresponding properties (ii), (iii), (iv), (vii) of Lemma 20, once we have shown that the probability measure $\Upsilon_N$ on $\mathbb{R}^2$ converges weakly to some limit $\Upsilon$. But (vi) of Lemma 20 implies tightness, so by passing to a further subsequence in $N$, we may assume $\Upsilon_N \to$ some $\Upsilon$. $\square$

4.6. *Completing the lower bound.* Now we will combine previous ingredients to show the following:

PROPOSITION 22. *For any $B'' > B$,*

$$\limsup_M M^{-2} c_{M,B''} \leq \gamma_*(\kappa, B) := \liminf_N N^{-2} \mathbb{E} \, \text{cost}_{(N)}(\mathbf{c}, B).$$

Lemma 11 shows

$$B \to \gamma_*(\kappa, B) \text{ is continuous nonincreasing on } \tfrac{1}{4} < B < \infty.$$

Now Propositions 17 and 22 combine to prove Theorem 2 and, hence, Theorem 1.

PROOF OF PROPOSITION 22. The proof mostly consists of combining Lemmas 19 and 21. Fix
$$M = KL; \qquad \delta > 0; \qquad B > \tfrac{1}{4}$$
and define
$$\varepsilon = 22B/L; \qquad B' = B + 16\delta L$$
and let $B'' \geq B' + \varepsilon$. First fix $D \subset \mathbb{R}^2$ and consider
$$Q^{[D]}(\cdot|\mathbf{c}) := \operatorname{tra}(\Xi_{\mathbf{c}}(\cdot \times D)),$$
where $\Xi_{\mathbf{c}}(\cdot,\cdot)$ is given in Lemma 21. Applying the linear map $Q \to Q_{\text{smooth}}$ to $Q^{[D]}(\cdot|\mathbf{c})$ gives transportation measures $Q^{[D]}_{\text{smooth}}(\cdot|\mathbf{c})$ across the extended $(M + 2K) \times (M + 2K)$ square such that
$$\operatorname{drift}(Q^{[D]}_{\text{smooth}}(\cdot|\mathbf{c})) - \operatorname{drift}(Q^{[D]}(\cdot|\mathbf{c})) = (x_{\mathbf{c},D}, y_{\mathbf{c},D}),$$
where, by (34), $|x_{\mathbf{c},D}| + |y_{\mathbf{c},D}| \leq \varepsilon$. If (say) $x_{\mathbf{c},D} > 0$ and $y_{\mathbf{c},D} < 0$, then in the notation of Lemma 12 (applied to $M + 2K$) the "adjusted" transportation measure
$$Q^{[D]}_{\text{adj}}(\cdot|\mathbf{c}) := Q^{[D]}_{\text{smooth}}(\cdot|\mathbf{c}) + x_{\mathbf{c},D}Q^{\leftarrow} + y_{\mathbf{c},D}Q^{\uparrow}$$
has
$$\operatorname{drift}(Q^{[D]}_{\text{adj}}(\cdot|\mathbf{c})) = \operatorname{drift}(Q^{[D]}(\cdot|\mathbf{c})). \tag{57}$$
The effect of the general $Q \to Q_{\text{smooth}}$ map on marginals is simply to average over blocks, and so the fact (55)
$$Q^{[D]} := \mathbb{E} Q^{[D]}(\cdot|\mathbf{c}) \in \mathcal{Q}_M$$
implies
$$Q^{[D]}_{\text{smooth}} := \mathbb{E} Q^{[D]}_{\text{smooth}}(\cdot|\mathbf{c}) \in \mathcal{Q}_{M+2K}.$$
So now
$$Q^{[D]}_{\text{adj}} := \mathbb{E} Q^{[D]}_{\text{adj}}(\cdot|\mathbf{c})$$
is of the form
$$Q^{[D]}_{\text{adj}} = Q^{[D]}_{\text{smooth}} + \varepsilon_{1,D}Q^{\rightarrow} + \varepsilon_{2,D}Q^{\leftarrow} + \varepsilon_{3,D}Q^{\uparrow} + \varepsilon_{4,D}Q^{\downarrow},$$
where $\sum_i \varepsilon_{i,D} \leq \varepsilon$. So clearly $Q^{[D]}_{\text{adj}} \in \mathcal{Q}_{M+2K}$. And
$$\operatorname{drift}(Q^{[D]}_{\text{adj}}) = \mathbb{E} \operatorname{drift}(Q^{[D]}_{\text{adj}}(\cdot|\mathbf{c}))$$
$$= \mathbb{E} \operatorname{drift}(Q^{[D]}(\cdot|\mathbf{c})) \qquad \text{by (57)}$$
$$= \operatorname{drift}(Q^{[D]})$$
$$= \int_D u\, \Upsilon(du) \qquad \text{by (56)}.$$

Writing $Q_{\text{adj}}$ for the $D = \mathbb{R}^2$ case of $Q_{\text{adj}}^{[D]}$, this density representation (cf. note below Lemma 21) shows $Q_{\text{adj}} \in \mathcal{Q}_{M+2K}^{\text{iso}}$. Apply Lemma 12 with $(Q_{\text{adj}}, Q_{\text{smooth}} := Q_{\text{smooth}}^{[\mathbb{R}^2]}, B'', B', M + 2K)$ in place of $(Q', Q, B', B, M)$:

$$(58) \quad \mathbb{E}\operatorname{cost}_{M+2K,B'}(Q_{\text{smooth}}, \mathbf{c}) \geq c_{M+2K,B''} - 4c^* B'' \varepsilon (M + 2K)^2.$$

Recalling the definition (33) of $m_B(Q)$,

$$m_{B'}(\mathbb{E} Q_{\text{smooth}}(\cdot|\mathbf{c})) = \mathbb{E}\operatorname{cost}_{M+2K,B'}(Q_{\text{smooth}}(\cdot|\mathbf{c}), \mathbf{c})$$
$$= \mathbb{E}\operatorname{cost}_{M+2K,B'}(Q_{\text{smooth}}, \mathbf{c}).$$

Now we apply Lemma 19 to the measures $\Xi_{\mathbf{c}}(\cdot \times \cdot)$ given by Lemma 21, with $\mathbf{f}_{\mathbf{c}}$ the associated flow-volume; conclusion (35) becomes

$$M^{-2}\mathbb{E}\operatorname{cost}_M(\mathbf{f}_{\mathbf{c}}, \mathbf{c}) \geq M^{-2}\mathbb{E}\operatorname{cost}_{M+2K,B'}(Q_{\text{smooth}}, \mathbf{c}) - 12c^* B^2/L$$
$$- 144B' c^* \delta L - 5c^* B^2 \ell(L, B/\delta) \exp(-M^2/L^2).$$

The left-hand side is upper bounded by $\gamma_*(\kappa, B)$ by (54). So combining this inequality with (58) gives

$$\gamma_*(\kappa, B) \geq M^{-2} c_{M+2K,B''} - 4c^* B'' \varepsilon \left(1 + \frac{4}{L}\right)^2$$
$$- 20c^* B^2/L - 256 B' c^* \delta L - 8c^* B^2 \ell(L, B/\delta) \exp(-M^2/L^2).$$

Now take $\delta = L^{-2}$ and then take $L = L_M \to \infty$ slowly with $M$. We get the conclusion of Proposition 22, except with the $\limsup_M$ taken through a subsequence $M_j$ with $M_{j+1}/M_j \to 1$. But then Lemma 15 identifies the $\limsup$ through that subsequence with the actual $\limsup$. □

REMARK. Without using Lemma 15, the argument gives the analog of Proposition 22 with $\liminf_M$ in place of $\limsup_M$, which in turn gives the (less aesthetic) analog of Theorem 2 with $\liminf_M$ in place of $\lim_M$.

**5. Discussion.** We start with comments tied closely to the statement and proof of Theorem 1 and then venture further afield.

1. We can apply the "method of bounded differences" argument (28) directly to the random cost $\operatorname{cost}_{(N)}(\mathbf{c}, B)$ in Theorem 1, and conclude

$$\mathbb{P}(|\operatorname{cost}_{(N)}(\mathbf{c}, B) - \mathbb{E}\operatorname{cost}_{(N)}(\mathbf{c}, B)| > \lambda) \leq 2 \exp\left(-\frac{\lambda^2}{4N^2(c^* B^2)^2}\right).$$

This of course implies the SLLN corresponding to Theorem 1. Various other "optimization over random data" problems have this character—that it is easy to show concentration about the mean, but not easy to show existence of a limit constant for the mean. Perhaps best known is random 3-SAT [8].

2. Our setup for Theorem 1 did not exclude the possibility $\mathbb{P}(c(e) = 0) > 0$. Recall [6, 12] that the critical value for bond percolation in $\mathbb{Z}^2$ equals $1/2$. It seems intuitively clear that

$$\text{if } \mathbb{P}(c(e) = 0) > \tfrac{1}{2}, \qquad \text{then } \gamma(\kappa, B) = 0;$$
$$\text{if } \mathbb{P}(c(e) = 0) < \tfrac{1}{2}, \qquad \text{then } \gamma(\kappa, B) > 0,$$

but we have not checked carefully.

3. The limit $\gamma(\kappa, B)$ is a priori nonincreasing in $B$. By considering random minimum-length paths (Lemma 6) and the bound $c^*$ on the support of $\kappa$, we find an upper bound

$$\gamma(\kappa, B) \leq c^*/8.$$

It is not hard to formalize the intuitive idea that as $B \downarrow 1/4$ the only feasible flows have $f(e) \approx 1/4$ for almost all $e$, and so

$$\lim_{B \downarrow 1/4} \gamma(\kappa, B) = \tfrac{1}{8} \mathbb{E} c(e).$$

4. Relaxing the requirement of bounded support, for what $\kappa$ do we still expect a finite limit $\gamma(\kappa, B)$? This question seems rather subtle. It is not obvious that the condition $\mathbb{E} c(e) < \infty$ is necessary, because we might be able to route flow to avoid edges $e$ with large $c(e)$. On the other hand, the condition $\mathbb{E} c^{1/4}(e) < \infty$ is necessary, otherwise [by considering the minimum of $c(e)$ over the four edges at a vertex] we have $\mathbb{E} \text{cost}_{(N)}(\mathbf{c}, B) = \infty$ for finite $N$.

5. Return to the case where $\kappa$ has bounded support. Define $\text{cost}_{(N)}(\mathbf{c}, \infty)$ as the "$B = \infty$" case where there is no bound on edge capacity. By monotonicity (of limit constants in $B$), the limit constant

$$\gamma(\kappa, \infty) := \lim_{B \to \infty} \gamma(\kappa, B)$$

exists. From Theorem 1 we get a bound

$$\limsup_N N^{-2} \mathbb{E} \text{cost}_{(N)}(\mathbf{c}, \infty) \leq \gamma(\kappa, \infty).$$

CONJECTURE 23. $\lim_N N^{-2} \mathbb{E} \text{cost}_{(N)}(\mathbf{c}, \infty) = \gamma(\kappa, \infty) = \lim_{M \to \infty} M^{-2} \times c_{M,\infty}$.

It seems plausible that one could just modify the arguments in this paper to directly study the $B = \infty$ case. Loosely speaking, in the $B = \infty$ case we expect that the optimal flow $f(e)$ across a typical edge $e$ has exponential tail

(59) $$\mathbb{P}(f(e) > x) \leq A_1 \exp(-x/A_2), \qquad 0 < x < \infty,$$

so that imposing a large bound $B$ has little effect.

6. It is tempting to speculate that the $Q \in \mathcal{Q}_M^{\mathrm{iso}}$ which attains the *inf* in the definition (10) should be some "geometrically natural" measure not depending on $\kappa$ or $B$. One possibility: take the intersection of Euclidean-isotropic lines with the continuous square $[0, M]^2$ to get a transportation measure on the boundary of this square; then discretize to get some isotropic transportation measure $\widetilde{Q}_M$ on the boundary of the discrete $M \times M$ square.

7. The construction in Section 3 gives flows of asymptotically optimal cost. From the use of the law of large numbers in the proof of Proposition 4, one can see that the routes of this flow (after scaling the $N \times N$ torus to the continuous torus $[0, 1)^2$) converge to straight line segments in the continuous torus. But we have not proved these are the *minimum-length* straight line segments; we have not excluded the possibility that some non-vanishing proportion of flow takes a long route around the torus instead of taking the shortest route. This possibility in taken into account in Lemma 21 (and the definition of $\mathcal{Q}_M^{\mathrm{iso}}$), where we use a general $\Upsilon$ rather than just the uniform distribution on $[-1/2, 1/2]^2$.

8. Given an environment **c** on the $N \times N$ torus, there is a calculus-type condition for a standardized global flow **f** to be a *local* minimum of the function $\mathbf{f} \to \mathrm{cost}_{(N)}(\mathbf{f}, \mathbf{c})$. This condition (*Waldrop equilibrium* [11]) is: for each source-destination pair $(x, y)$, the flow from $x$ to $y$ can only use minimum-weight paths, where here the *weight* of an edge $e$ is $c(e)f(e)$. But we do not know how to exploit this condition in studying the optimal flow.

9. If instead of the quadratic costs (1) we used linear costs $\sum_e c(e)f(e)$, then the optimal flow simply chooses the minimum-cost path [where cost of a path equals sum of edge-costs $c(e)$] between each source and each destination; taking $B = \infty$, there is no interaction between flows. This case is just first passage percolation, and mean costs in our "flow" setting are easily related to the time constant in first passage percolation. Questions concerning flow volume across edges have apparently not been studied rigorously, but one expects power-law tail behavior instead of (59).

10. There is a huge literature on algorithms for flows in networks, illustrated by the monograph [1]. But the kind of multicommodity flow problems that we study are typically NP-hard as algorithmic problems. There is a large body of theoretical work going back to Leighton–Rao [15] relating multicommodity flow to other network problems, but this typically gives upper and lower bounds differing by log(number of vertices) factors. We should emphasize that the basic max-flow min-cut theorem does *not* apply in our multicommodity flow setting.

11. The kind of feasible flow problems mentioned in the introduction have been studied in the *unidirectional* case where one seeks to maximize volume of flow from somewhere on the bottom side of a square to somewhere on the top side. See [7, 9, 13, 21] for results and connections with percolation and first-passage percolation. In this setting one *can* apply the max-flow min-cut theorem. Otherwise there seems no literature very closely related to Theorem 1, though of course the general idea "consider large blocks" is pervasive throughout the study of spatial stochastic processes.

12. A more abstract view of our argument goes as follows. Subsequential weak limits as $N \to \infty$ of optimal flows across a typical $M \times M$ square define flows across $M \times M$ squares in the lattice $\mathbb{Z}^2$. These are consistent as $M$ increases and so define a flow on $\mathbb{Z}^2$ "from the infinite boundary to the infinite boundary." The constant $\gamma(\kappa, B)$ defined in Theorem 2 as a $M \to \infty$ limit could more abstractly be defined in terms of optimal costs in an appropriate problem involving such infinite-distance flows on $\mathbb{Z}^2$. See [3] for the details of such an approach in a different $\mathbb{R}^2$ setting (traveling salesman problem through some specified proportion of random points).

13. The abstract view above does not seem to help in the $d$-dimensional setting, but if instead the random network is "locally tree-like" (as with many models of random graph), it can be related to the *cavity method* [16] of statistical physics. That is, the $n \to \infty$ limit problem is an optimization problem on an infinite *tree* which can be tackled by setting up recursive equations exploiting the tree structure. Making such arguments rigorous is a major challenge; see [2] for nonrigorous methodology applied to optimal flow problems and leading to explicit numerical results.

14. What kind of flow problems can be handled by the methodology of this paper, as outlined in Section 2.1? Intuitively, what seems important is the following:

- We are studying an optimal global flow which minimizes a global cost function defined as a sum of local cost functions.
- Constraints are local.
- There is a stationary random environment determining costs and constraints.

In such a problem one can seek to modify Theorem 2 by replacing the specific definition of $c_{M,B}$ with its analog for different model assumptions. But we need the solution to have the property that flow volumes across different edges have the same order of magnitude, eliminating the "linear" case in comment 9 above.

**Acknowledgments.** I thank an anonymous referee and Associate Editor for detailed comments and encouragement to improve the outline of the proof, and Sourav Chatterjee for helpful discussion.

Department of Statistics
367 Evans Hall # 3860
UC Berkeley
Berkeley, California 94720
E-mail: aldous@stat.berkeley.edu
URL: www.stat.berkeley.edu/users/aldous